\documentclass{amsproc}

\usepackage{graphicx}

\newtheorem{theorem}{Theorem}[section]
\newtheorem{lemma}[theorem]{Lemma}
\newtheorem{cor}[theorem]{Corollary}
\newtheorem{prop}[theorem]{Proposition}

\theoremstyle{definition}

\theoremstyle{remark}

\numberwithin{equation}{section}

\begin{document}

\title[Alternating Permutations]{A Survey of Alternating
  Permutations} 


\author{Richard P. Stanley}
\address{Department of Mathematics, M.I.T., Cambridge, MA 02139}
\curraddr{}
\email{rstan@math.mit.edu}
\thanks{This material is based upon work supported by the
  National Science Foundation under Grant No.~0604423. Any opinions,
  findings and conclusions or recommendations expressed in this
  material are those of the author and do not necessarily reflect
  those of the National Science Foundation.}


\date{March 18, 2009}

\begin{abstract} 
A permutation $a_1a_2\cdots a_n$ of $1,2,\dots,n$ is
\emph{alternating} if $a_1>a_2<a_3>a_4<\cdots$. We survey some aspects
of the theory of alternating permutations, beginning with the famous
result of Andr\'e that if $E_n$ is the number of alternating
permutations of $1,2,\dots,n$, then $\sum_{n\geq 0}E_n\frac{x^n}{n!} =
\sec x +\tan x$. Topics include refinements and $q$-analogues of
$E_n$, various occurrences of $E_n$ in mathematics, longest
alternating subsequences of permutations, umbral enumeration of
special classes of alternating permutations, and the connection
between alternating permutations and the $cd$-index of the symmetric
group. 
\end{abstract}

\maketitle

\begin{centering}
{\small  Dedicated to Reza Khosrovshahi on the occasion of his
     70th birthday}\\ 
\end{centering}


\bibliographystyle{amsplain}

\newcommand{\beq}{\begin{equation}} 
\newcommand{\eeq}{\end{equation}}
\newcommand{\zz}{\mathbb{Z}}
\newcommand{\rr}{\mathbb{R}}
\newcommand{\cc}{\mathbb{C}}
\newcommand{\cp}{\mathcal{P}}
\newcommand{\cac}{\mathcal{C}}
\newcommand{\sn}{\mathfrak{S}_n} 
\newcommand{\fs}{\mathfrak{S}}
\newcommand{\fo}{\mathfrak{o}}
\newcommand{\st}{\,:\,} 
\newcommand{\hy}{\!-\!}
\newcommand{\as}{\mathrm{as}}
\newcommand{\alt}{\mathrm{Alt}}
\newcommand{\ralt}{\mathrm{Ralt}}
\newcommand{\maj}{\mathrm{maj}}
\newcommand{\inv}{\mathrm{inv}}
\newcommand{\is}{\mathrm{is}}
\newcommand{\bea}{\begin{eqnarray}}
\newcommand{\eea}{\end{eqnarray}}
\newcommand{\be}{\begin{enumerate}}
\newcommand{\ee}{\end{enumerate}}
\newcommand{\beas}{\begin{eqnarray*}}
\newcommand{\eeas}{\end{eqnarray*}}
\newcommand{\hz}{\hat{0}}
\newcommand{\ho}{\hat{1}}

\thispagestyle{empty}

\vskip 20pt

\section{Basic enumerative properties.} \label{sec1} Let $\sn$ denote
the symmetric group of all permutations of $[n]:=\{1,2,\dots,n\}$. A
permutation $w=a_1a_2\cdots a_n\in\sn$ is called \emph{alternating} if
$a_1>a_2<a_3>a_4<\cdots$. In other words, $a_i<a_{i+1}$ for $i$ even,
and $a_i>a_{i+1}$ for $i$ odd. Similarly $w$ is \emph{reverse
  alternating} if $a_1<a_2>a_3<a_4>\cdots$. (Some authors reverse
these definitions.) Let $E_n$ denote the
number of alternating permutations in $\sn$. (Set $E_0=1$.) For
instance, $E_4=5$, corresponding to the permutations 2143, 3142, 3241,
4132, and 4231.  The number $E_n$ is called an \emph{Euler number}
because Euler considered the numbers $E_{2n+1}$, though not with the
combinatorial definition just given. (Rather, Euler defined them via
equation~\eqref{eq:tangf} below.) The Euler numbers are not to be
confused with the \emph{Eulerian numbers}, which count permutations by
number of descent. The involution 
  \beq a_1 a_2 \cdots a_n \mapsto n+1-a_1,n+1-a_2,\cdots,n+1-a_n  
    \label{eq:inv} \eeq
on $\sn$ shows that $E_n$ is also the number of \emph{reverse}
alternating permutations in $\sn$. We write $\alt_n$ (respectively,
$\ralt_n$) for the set of alternating (respectively, reverse
alternating) permutations $w\in\sn$.

The subject of alternating permutations and Euler numbers has become
so vast that it is impossible to give a comprehensive survey. We will
confine ourselves to some highlights and to some special topics that
we find especially interesting.

The fundamental enumerative property of alternating permutations is
due to Desir\'e Andr\'e \cite{andre} in 1879.  (Note however that
Ginsburg \cite{ginsburg} asserts without giving a reference that Binet
was aware before Andr\'e that the coefficients of $\sec x$ count
alternating permutations.)

\begin{theorem} \label{thm:gf}
We have
  \beas \sum_{n\geq 0}E_n\frac{x^n}{n!} & = & \sec x + \tan x\\
    & = & 1+x+\frac{x^2}{2!}+2\frac{x^3}{3!}+5\frac{x^4}{4!}
    +16\frac{x^5}{5!}+61\frac{x^6}{6!}+ 272\frac{x^7}{7!}
    +1385\frac{x^8}{8!}+\cdots. \eeas
\end{theorem}

Note that $\sec x$ is an even function (i.e, $\sec(-x)=\sec x$),
while $\tan x$ is odd ($\tan(-x)=-\tan x$). It follows from
Theorem~\ref{thm:gf} that
  \bea \sum_{n\geq 0}E_{2n}\frac{x^{2n}}{(2n)!} & = & \sec x
     \label{eq:secgf}\\
     \sum_{n\geq 0}E_{2n+1}\frac{x^{2n+1}}{(2n+1)!} & = & \tan x. 
    \label{eq:tangf}\eea
For this reason $E_{2n}$ is sometimes called a \emph{secant number}
and $E_{2n+1}$ a \emph{tangent number}.

We will sketch three proofs of Theorem~\ref{thm:gf}.

\emph{First proof.} Let $0\leq k\leq n$. Choose a $k$-subset $S$ of
$[n]=\{1,2,\dots,n\}$ in 
$\binom{n}{k}$ ways, and set $\bar{S}=[n]-S$. Choose a reverse
alternating permutation $u$ of $S$ in $E_k$ ways, and choose a reverse
alternating permutation $v$ of $\bar{S}$ in $E_{n-k}$ ways. Let $w$ be
the concatenation $u^r,n+1,v$, where $u^r$ denotes the reverse of $u$
(i.e., if $u=u_1\cdots u_k$ then $u^r=u_k\cdots u_1$). When $n\geq 2$,
we obtain in this way every alternating and every reverse alternating
permutation $w$ exactly once. Since there is a bijection between
alternating and reverse alternating permutations of any finite
(ordered) set, the number of $w$ obtained is $2E_{n+1}$. Hence
  \beq 2E_{n+1} = \sum_{k=0}^n \binom nk E_k E_{n-k},\ \ n\geq 1. 
    \label{eq:eulerrec} \eeq
Set $y= \sum_{n\geq 0}E_n x^n/n!$. Taking into account the initial
conditions $E_0=E_1=1$, equation~\eqref{eq:eulerrec} becomes the
differential equation
  $$ 2y' = y^2 +1,\ \ y(0)=1. $$
The unique solution is $y=\sec x+\tan x$.
\qed

\textsc{Note.} The clever counting of both alternating and reverse
alternating permutations in the proof of
Theorem~\ref{thm:gf} can be avoided at the cost of a little
elegance. Namely, by considering the position of 1 in an alternating
permutation $w$, we obtain the recurrence
  $$ E_{n+1}=\sum_{\substack{1\leq j\leq n\\ j\ \mathrm{odd}}}
    \binom nj E_jE_{n-j},\ \ n\geq 1. $$
This recurrence leads to a system of differential equations for the
power series $\sum_{n\geq 0}E_{2n}x^{2n}/(2n)!$ and $\sum_{n\geq
  0}E_{2n+1}x^{2n+1}/(2n+1)!$.

\emph{Second proof.} For simplicity we consider only $E_{2n}$. A
similar, though somewhat more complicated, argument works for
$E_{2n+1}$. We want to show that
     $$ \left(\sum_{n\geq 0} E_{2n}\frac{x^n}{n!}\right)\left(
       1-\frac{x^2}{2!}+\frac{x^4}{4!}-\cdots\right) =1. $$
   Equating coefficients of $x^{2n}/(2n)!$ on both sides gives
   \beq E_{2n} = \binom n2 E_{2n-2}-\binom n4 E_{2n-4}+
     \binom n6 E_{2n-6}-\cdots. \label{eq:e2nie} \eeq
 Let $S_k$ be the set of
   permutations $w=a_1 a_2\cdots a_{2n}\in\fs_{2n}$ satisfying
    $$ a_1>a_2<a_3>a_4<\cdots >a_{2n-2k},\ \ a_{2n-2k+1}>
       a_{2n-2k+2}>\cdots>a_{2n}, $$
   and let $T_k$ be those permutations in $S_k$ that also satisfy 
   $a_{2n-2k}>a_{2n-2k+1}$. Hence $S_1-T_1$ consists of all
   alternating permutations in $\sn$. Moreover,
   $T_i=S_{i+1}-T_{i+1}$. Hence
    $$ E_n = \#(S_1-T_1)=\#S_1-\#(S_2-T_2)= \cdots = \#S_1-\#S_2+
     \#S_3-\cdots. $$
   A permutation in $S_k$ is obtained by choosing $a_{2n-2k+1}$,
   $a_{2n-2k+2},\dots$, $a_{2n}$ in $\binom{2n}{2k}$ ways and then
   $a_1, a_2,\dots,a_{2n-2k}$ in $E_{2(n-k)}$ ways. Hence $\#S_k =
   \binom{2n}{2k}E_{2(n-k)}$, and the proof follows. \qed

\emph{Third proof.} Our third proof gives a more general formula that
makes it more obvious why $\sec x$ appears in Theorem~\ref{thm:gf}.
Our second proof can also be extended to yield equation~\eqref{eq:fkx}
below.  A slightly more complicated argument, omitted here, explains
the term $\tan x$. For some generalizations, see \cite{g-j} and
\cite[Exer.~3.80]{ec1}.

Let $k,n\geq 1$, and let $f_k(n)$ denote the number of permutations
$a_1 a_2\cdots a_{kn}\in \fs_{k n}$ satisfying $a_i>a_{i+1}$ if and
only if $k|i$. Set 
  $$ F_k(x) = \sum_{n\geq 0} f_k(n)\frac{x^{kn}}{(kn)!}. $$
We claim that
  \beq F_k(x) = \frac{1}{\sum_{n\geq 0}(-1)^n\frac{x^{kn}}{(kn)!}}. 
     \label{eq:fkx} \eeq
To prove equation~\eqref{eq:fkx}, note that the number of permutations
$b_1 b_2\cdots b_{kn}\in \fs_{kn}$ such that $b_i>b_{i+1}$ only if
(but not necessarily if) $k|i$ is easily seen to be the multinomial
coefficient $\binom{kn}{k,k,\dots,k} = (kn)!/k!^n$. A straightforward
inclusion-exclusion argument then gives 
   \beq f_k(n) = \sum_{j=1}^n \sum_{\substack{i_1+\cdots+i_j=n\\ i_r>0}}
   (-1)^{n-j} \binom{kn}{i_1 k,\dots,i_j k}. \label{eq:incex} \eeq
Comparing with the expansion
   $$  \frac{1}{\sum_{n\geq 0}(-1)^n\frac{x^{kn}}{(kn)!}} =
    \sum_{j\geq 0}\left(\sum_{n\geq 1}(-1)^{n-1} \frac{x^{kn}}{(kn)!}
    \right)^j $$
completes the proof. 
\qed

The inclusion-exclusion formula \eqref{eq:incex} can be regarded as
the expansion of a determinant, giving the determinantal formula
  \beq f_k(n)=(kn)!\det\left[1/(k(j-i+1))!\right]_{i,j=1}^n, 
       \label{eq:fkndet} \eeq
where we set $1/(-m)!=0$ for $m>0$. (See \cite[pp.~69--71]{ec1}.) The
case $k=2$ gives a formula for $E_{2n}$. Similarly there holds
  \beq E_{2n-1} =(2n-1)!\det\left[ 1/c_{ij}!\right]_{i,j=1}^n, 
        \label{eq:oddedet} \eeq
where $c_{1,j}=2j-1$ and $c_{ij}=2(j-i+1)$ for $2\leq i\leq n$. 

Because $f(z)=\sec z+\tan z$ is a well-behaved function of the complex 
variable $z$ (in particular, it is meromorphic, and all poles are
simple), it is routine to derive a precise asymptotic estimate of
$E_n$. The smallest pole of $f(z)$ is at $z=\pi/2$,
with residue $-2$, and the next smallest pole is at
$z=-3\pi/2$. Hence
   \beq \frac{E_n}{n!} = \frac{4}{\pi}\left(\frac{2}{\pi}\right)^n
                 +O\left(\left(\frac{2}{3\pi}\right)^n\right). 
      \label{eq:eulasymp} \eeq
In fact, the poles of $f(z)$ are precisely $z=(-1)^n(2n+1)\pi/2$,
$n\geq 0$, all with residue $-2$, leading to the \emph{convergent}
asymptotic series
  $$ \frac{E_n}{n!}=2\left(\frac{2}{\pi}\right)^{n+1}
        \sum_{k\geq 0}(-1)^{k(n+1)} 
        \frac{1}{(2k+1)^{n+1}}. $$
This formula for $n$ odd is equivalent to the well-known evaluation of
$\zeta(n+1):=\sum_{k\geq 1}k^{-(n+1)}$. For further information on
asymptotic expansions see \cite{f-s}, where in particular Euler
numbers are treated in Example~IV.35.

\section{Refinements of Euler numbers} \label{sec:refine}

For the purposes of this paper, a \emph{refinement} of the Euler
number $E_n$ is a sequence $a_0,a_1,\dots $of nonnegative integers
summing to $E_n$. Often we encode the sequence by a polynomial
$P(q)=\sum a_iq^i$, so $P(1)=E_n$. This polynomial is sometimes called
a \emph{$q$-analogue} of $E_n$, especially if it involves such
``$q$-objects'' as the finite field $\mathbb{F}_q$ or expressions such
as $[n]!:=(1-q)(1-q^2)\cdots (1-q^n)$ or $\sum a_n x^n/[n]!$. We
briefly discuss two refinements of $E_n$.

The first refinement provides an elegant scheme for computing
the Euler numbers. Let $E_{n,k}$ denote the number of alternating
permutations of $[n+1]$ with first term $k+1$. For instance,
$E_{n,n}=E_n$. It is easy to verify the recurrence
  $$ E_{0,0}=1,\ \ E_{n,0}=0\ (n\geq 1),\ \  E_{n+1,k+1}
    =E_{n+1,k}+E_{n,n-k}\ (n\geq k\geq 0). $$
   Note that if we place the $E_{n,k}$'s in the triangular array
   \beq \begin{array}{ccccccccc} & & & & E_{00}\\ & & & E_{10} &
     \rightarrow & 
     E_{11}\\ & & E_{22} & \leftarrow & E_{21} & \leftarrow & E_{20}\\
    & E_{30} & \rightarrow & E_{31} & \rightarrow & E_{32} & \rightarrow
    & E_{33}\\ E_{44} & \leftarrow & E_{43} & \leftarrow & E_{42} &
    \leftarrow & E_{41} & \leftarrow & E_{40}\\
    & & & & \cdots \end{array} \label{eq:boustri} \eeq
   and read the entries in the direction of the arrows from
   top-to-bottom (the so-called \emph{boustrophedon} 
  or \emph{ox-plowing} order), then the first
   number read in each row is 0, and each subsequent entry is the sum
   of the previous entry and the entry above in the previous row. The
   first seven rows of the array are as follows:
  $$ \begin{array}{ccccccccccccc} & & & & & & 1\\ & & & & & 0 &
     \rightarrow & 
     1\\ & & & & 1 & \leftarrow & 1 & \leftarrow & 0\\
    & & & 0 & \rightarrow & 1 & \rightarrow & 2 & \rightarrow
    & 2\\ & & 5 & \leftarrow & 5 & \leftarrow & 4 &
    \leftarrow & 2 & \leftarrow & 0\\ & 0 & \rightarrow & 5
     & \rightarrow & 10 & \rightarrow & 14 & \rightarrow & 
     16 & \rightarrow & 16\\ 61 & \leftarrow & 61 & \leftarrow & 
     56 & \leftarrow & 46 & \leftarrow & 32 & \leftarrow & 
     16 & \leftarrow & 0.\\
    & & & & & & \cdots \end{array} $$

We can obtain a generating function for the number $E_{n,k}$ as
follows.  Define
   $$ [m,n] = \left\{ \begin{array}{rl} m, & m+n\ \mathrm{odd}\\
      n, & m+n\ \mathrm{even}. \end{array} \right. $$
Then
  \beq \sum_{m\geq 0}\sum_{n\geq 0} E_{m+n,[m,n]} \frac{x^m}{m!}
     \frac{y^n}{n!} = \frac{\cos x+\sin x}{\cos(x+y)}. 
     \label{eq:entgf} \eeq
For a proof see Graham, Knuth, and Patashnik
\cite[Exer.~6.75]{g-k-p}.  The numbers $E_{n,k}$ are called
\emph{Entringer numbers}, after R. C. Entringer \cite{ent}. The
triangular array \eqref{eq:boustri} is due to L. Seidel \cite{seidel}
(who used the word ``boustrophedon'' to describe the triangle).  It
was rediscovered by Kempner \cite{kempner}, Entringer \cite{ent}
and Arnold \cite{arnold}.  For further information and
references, see J. Millar, N. J. A.  Sloane, and N. E. Young
\cite{m-s-y}. A more recent reference is R. Ehrenborg and S. Mahajan
\cite[{\S}2]{eh-ma}.

Our second refinement $E_n(q)$ of Euler numbers is a natural
$q$-analogue. An \emph{inversion} of a permutation $w=a_1\cdots
a_n\in\sn$ is a pair $(i,j)$ such that $i<j$ and $a_i>a_j$. Let
inv$(w)$ denote the number of inversions of $w$, and define
  \beq E_n(q) = \sum_{w\in\ralt_n} q^{\inv(w)}. \label{eq:enqdef} \eeq  
For instance, we have inv$(1324)=1$, inv$(1423)=2$,
inv$(2314)=2$, inv$(2413)=3$, and inv$(3412)=4$, so
$E_n(q)=q+2q^2+q^3+q^4$. Similarly define
  $$ E^*_n(q) = \sum_{w\in\alt_n } q^{\inv(w)}. $$
For instance, we have inv$(2143)=2$, inv$(3142)=3$, inv$(3241)=4$,
inv$(4132)=4$, and inv$(4231)=5$, so $E_4(q)=q^2+q^3+2q^4+q^5$. 
Note that 
   \beq E^*_n(q) = q^{\binom n2}E_n(1/q), \label{eq:enq} \eeq
an immediate consequence of the involution~\eqref{eq:inv}. The
polynomials $E_n(q)$ for $2\leq n\leq 6$ are given by
    \beas E_2(q) & = & 1\\
          E_3(q) & = & q+q^2\\
          E_4(q) & = & q^2+q^3+2q^4+q^5\\
          E_5(q) & = & q^2+2q^3+3q^4+4q^5+3q^6+2q^7+q^8\\
          E_6(q) & = & q^3+2q^4+5q^5+7q^6+9q^7+10q^8+10q^9+8q^{10}
                       +5q^{11}+2q^{12}+q^{13}. \eeas

The main combinatorial result on the ``$q$-Euler polynomials''
$E_n(q)$ and $E^*_n(q)$, known to Sch\"utzenberger in 1975 (according
to D.~Foata) and rediscovered by Stanley \cite{rs:binpos} in the
context of binomial posets, is the following. Define the $q$-cosine
and $q$-sine functions by 
   \beas \cos_q(x) & = & \sum_{n\geq 0}(-1)^n
   \frac{x^{2n}}{[2n]!}\\
        \sin_q(x) & = & \sum_{n\geq 0}(-1)^n
        \frac{x^{2n+1}}{[2n+1]!}, \eeas
where $[m]!=(1-q)(1-q^2)\cdots(1-q^m)$. 
Note that these series become $\cos x$ and $\sin x$ after
substituting $(1-q)x$ for $x$ and letting $q\to 1$. Similarly define
the variants
    \beas \cos^*_q(x) & = & \sum_{n\geq 0}(-1)^n q^{\binom{2n}{2}}
   \frac{x^{2n}}{[2n]!}\\
        \sin^*_q(x) & = & \sum_{n\geq 0}(-1)^n q^{\binom{2n+1}{2}}
        \frac{x^{2n+1}}{[2n+1]!}. \eeas
Note that 
   \beq \cos^*_q(x) = \cos_{1/q}(-x/q), \label{eq:qcos} \eeq
and similarly for $\sin^*_q(x)$.

\begin{theorem} \label{thm:qeuler}
We have
  \beas \sum_{n\geq 0}E_n(q)\frac{x^n}{[n]!} & = & 
           \frac{1}{\cos_q(x)}+\frac{\sin_q(x)}{\cos_q(x)}\\
  \sum_{n\geq 0}E^*_n(q)\frac{x^n}{[n]!} & = & 
        \frac{1}{\cos^*_q(x)}+\frac{\sin^*_q(x)}{\cos^*_q(x)}   
     \eeas
\end{theorem}

One way to prove Theorem~\ref{thm:qeuler} is by a straightforward
generalization of the third proof of Theorem~\ref{thm:gf}.  A more
conceptual explanation for this result and some generalizations based
on binomial posets appear in \cite{rs:binpos}. Note that the result
for $E^*_n(q)$ is an immediate consequence of that for $E_n(q)$ and
equations~\eqref{eq:enq} and \eqref{eq:qcos}. The reasoning used to
deduce the determinantal formulas \eqref{eq:fkndet} and
\eqref{eq:oddedet} from equation~\eqref{eq:incex} can also be
straightforwardly generalized to give
  \bea E_{2n}(q) & = &  [2n-1]!\det\left[1/[2(j-i+1)]!
    \right]_{i,j=1}^n \label{eq:enqdete}\\
     E_{2n-1}(q) & = & [2n-1]!\det\left[ 1/[c_{ij}]!
        \right]_{i,j=1}^n, \label{eq:enqdeto} \eea
where $c_{1,j}=2j-1$ and $c_{ij}=2(j-i+1)$ for $2\leq i\leq n$. 
 
The polynomials $E_n(q)$ also appear in a natural way in the theory of
symmetric functions. Assuming familiarlity with this theory
\cite{mac}\cite[Ch.~7]{ec2}, let $\tau_n$ be as in
Section~\ref{subsec:zzt} of this paper, and let $s_{\tau_n}$ be the
corresponding skew Schur function. Then
  $$ E_n(q)=[n]!s_{\tau_n}(1,q,q^2,\dots). $$
It follows from standard properties of Schur functions (see
\cite[Prop.~7.19.11]{ec2}) that $E_n(q)$ has the additional
combinatorial interpretation
  \beq E_n(q) = \sum_{w\in\alt_n} q^{\maj(w^{-1})}, \label{eq:enqmaj}
  \eeq 
where 
   $$ \maj(w) = \sum_{\substack{1\leq i\leq n-1\\ a_i>a_{i+1}}}i, $$
the \emph{major index} of $w=a_1\cdots a_n$. The equivalence of
equations~\eqref{eq:enqdef} and \eqref{eq:enqmaj} is also a
consequence of well-known properties of inv and maj (e.g.,
\cite[Thm.~1.4.8]{ec1-2}). In the context of symmetric functions, the
determinantal formulas \eqref{eq:enqdete} and \eqref{eq:enqdeto} are
consequences of the Jacobi-Trudi identity \cite[{\S}7.16]{ec2} for the
skew Schur function $s_{\tau_n}$. 

A number of other $q$-analogues of $E_n$ (sometimes just for $n$ odd)
have been proposed; see
\cite{fo-ha}\cite{fulmek}\cite{h-r-z}\cite{prod1}\cite{prod2} for
further information.

\section{Other occurrences of Euler numbers} \label{sec2}
There are numerous occurrences of Euler numbers not directly related
to alternating permuations. We discuss a few of these here. For more
information on this topic, see the treatise \cite{viennot2} of
Viennot. 

\subsection{Complete increasing binary trees}
\begin{figure}
\centering
\centerline{\includegraphics[width=8cm]{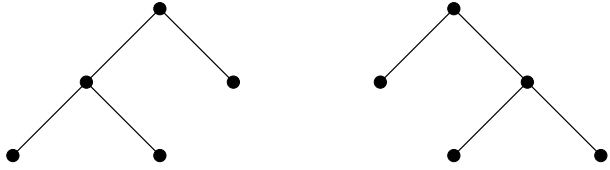}}
\caption{The two complete binary trees with five vertices}
\label{fig:ibt}
\end{figure}

A \emph{(plane) binary tree} on the vertex set $[n]$ is defined
recursively by having a root vertex $v$ and a left and right subtree
of $v$ which are themselves binary trees or are empty. A binary tree
is \emph{complete} if every vertex either has two children or is an
endpoint. A binary tree on the vertex set $[n]$ is \emph{increasing}
if every path from the root is increasing. Figure~\ref{fig:ibt} shows
the two complete binary trees with five vertices. Each one has
eight increasing labelings, so there are 16 complete increasing binary
trees on the vertex set $[5]$.

\begin{theorem} \label{prop:ibt}
The number of complete increasing binary trees on $[2m+1]$ is the
Euler number $E_{2m+1}$. (There is a similar but more complicated
statement for the vertex set $[2m]$ which we do not give here.)
\end{theorem}

\proof
Given a sequence $u=a_1 a_2\cdots a_n$ of distinct integers, define a
labelled binary tree $T_u$ as follows. Let
$a_i=\min\{a_1,\dots,a_n\}$. Let $a_i$ be the root of $T_w$, and
recursively define the left subtree of the root to be $T_{a_1\cdots
  a_{i-1}}$, and the right subtree of the root to be $T_{a_{i+1}\cdots
  a_n}$. It is not hard to check that the map $w\mapsto T_w$ is a
bijection from alternating permutations $w\in\fs_{2m+1}$ to complete
increasing binary trees on $[2m+1]$.
\qed

\subsection{Flip equivalence} \label{subsec:flip}
\begin{figure}
\centering
\centerline{\includegraphics[width=12cm]{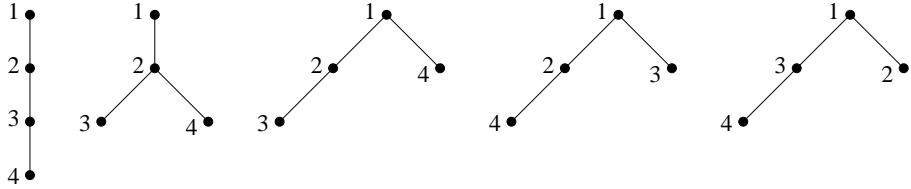}}
\caption{The five increasing 1-2 trees with four vertices}
\label{fig:5ibt}
\end{figure}

The Euler numbers are related to increasing binary trees in another
way. The \emph{flip} of a binary tree at a vertex $v$ is the binary
tree obtained by interchanging the left and right subtrees of $v$.
Define two increasing binary trees $T$ and $T'$ on the vertex set
$[n]$ to be \emph{flip equivalent} if $T'$ can be obtained from $T$ by
a sequence of flips. Clearly this definition of flip equivalence is an
equivalence relation. The equivalence classes are in an obvious
bijection with increasing \emph{1-2 trees} on the vertex set $[n]$, that
is, increasing (rooted) trees so that every non-endpoint vertex has
one or two children.  (These are not plane trees, i.e., the order in
which we write the children of a vertex is irrelevant.)
Figure~\ref{fig:5ibt} shows the five increasing 
1-2 trees on four vertices, so $f(4)=5$. 

\begin{theorem} \label{prop:fostr}
We have $f(n)=E_n$ (an Euler number).
\end{theorem}

\proof
Perhaps the most straightforward proof is by generating functions. (A
combinatorial proof was first given by Donaghey \cite{don}.) Let
  $$ y=\sum_{n\geq 1}f(n)\frac{x^n}{n!} =
  x+\frac{x^2}{2}+2\frac{x^3}{6} +\cdots. $$
Then $y'=\sum_{n\geq 0}f(n+1)x^n/n!$. Every increasing 1-2 tree with
$n+1$ vertices is either (a) a single vertex ($n=0$), (b) has one
subtree of the root which is an increasing binary tree with $n$
vertices, or (c) has two subtrees of the root, each of which is an
increasing binary tree, with $n$ vertices in all. The order of the two
subtrees is irrelevant. From this
observation we obtain the differential equation $y'=1+y+\frac 12y^2$,
$y(0)=0$. The unique solution is $y=\sec x+\tan x-1$, and the proof
follows from Theorem~\ref{thm:gf}.
\qed

Much additional information concerning the connection between
alternating permutations and increasing trees appears in a paper of
Kuznetsov, Pak, and Postnikov \cite{k-p-p}.

\textsc{Algebraic note.} Let $\mathcal{T}_n$ be the set of all
increasing binary trees with vertex set $[n]$. For $T\in\mathcal{T}_n$
and $1\leq i\leq n$, let $\omega_iT$ be the flip of $T$ at vertex $i$.
Then clearly the $\omega_i$'s generate a group isomorphic to
$(\zz/2\zz)^n$ acting on $\mathcal{T}_n$, and the orbits of this
action are the flip equivalence classes. For further details see Foata
\cite{foata} and Foata-Strehl \cite{fo-st}.

\subsection{Ballot sequences}
A \emph{ballot sequence} of length $2n$ is a sequence
$b_1,b_2,\dots$, $b_{2n}$ of $n$ $1$'s and $n$ $-1$'s for which all
partial sums are nonnegative. The number of ballot sequences of
length $2n$ is the Catalan number $C_n$
\cite[Cor.~6.2.3(ii)]{ec2}. Given a ballot sequence
$b=b_1,b_2,\dots,b_{2n}$, define 
  \beas \omega(b) & = & \prod_{i\st b_i=1}(b_1+b_2+\cdots +b_i)\\
  \omega^*(b) & = & \prod_{i\st b_i=1}(b_1+b_2+\cdots +b_i+1). \eeas
For instance, if $b=(1,1,-1,1,-1,-1,1,-1)$ then $\omega(b)=1\cdot
2\cdot 2\cdot 1=4$ and $\omega^*(b)=2\cdot 3\cdot 3\cdot 2=
36$. Let $B(n)$ denote the set of all ballot sequences of length 
$2n$. It follows from a bijection of Fran\c{c}on and Viennot
\cite{fr-vi} between binary trees and weighted paths (also explained
in \cite[{\S}5.2]{go-ja}) that 
  \bea \sum_{b\in B(n)} \omega(b)^2 & = &  E_{2n}. \label{eq:dyck}\\
      \sum_{b\in B(n)} \omega(b)\omega^*(b) & = &  E_{2n+1}.    
     \nonumber \eea

Equation~\eqref{eq:dyck} suggests the following question: is there a
finite-dimensional algebra $\mathcal{A}_n$ (say over $\cc$) with a
``natural'' combinatorial definition whose irreducible representations
have dimension $\omega(b)$, $b\in B(n)$? It would then follow that
$\dim\mathcal{A}_n/\sqrt{\mathcal{A}_n}= E_{2n}$, where
$\sqrt{\mathcal{A}_n}$ denotes the radical of $\mathcal{A}_n$.

\subsection{Simsun permutations} \label{ssec:simsun}
 Define a \emph{simsun permutation} 
  to be a permutation $w=a_1\cdots a_n\in\sn$ such that for all $1\leq
  k\leq n$, the subword of $w$ consisting of $1,2,\dots,k$ (in the
  order they appear in $w$) does not have three consecutive decreasing
  elements.  For instance, $w=425163$ is not simsun since if we remove
  5 and 6 from $w$ we obtain $4213$, which contains the three
  consecutive decreasing elements 421. Simsun permutations were named
  after Rodica Simion and Sheila Sundaram and were first described in
  print by Sundaram \cite[{\S}3]{sund}. They are a variant of a class
  of permutations due to Foata and Sch\"utzenberger \cite{fo-sc} known
  as \emph{Andr\'e permutations}. We have chosen here to deal only
  with simsun permutations because their definition is a little
  simpler than that of Andr\'e permutations. Simion and Sundaram prove
  in their paper the following basic result on simsum permutations.

\begin{theorem} \label{thm:simsun}
The number of simsun permutations in $\sn$ is the Euler number $E_n$. 
\end{theorem}

\textsc{Proof} (sketch).  Let $f_k(n)$ be the number of
   simsun permutations $w=a_1\cdots a_n\in\sn$ with $k$ descents,
   i.e., $k$ values of $i$ for which $a_i>a_{i+1}$. By inserting $n+1$ 
   into a simsun permutation in $\sn$, we get the recurrence
  $$ f_k(n+1) = (n-2k+2)f_{k-1}(n)+(k+1)f_k(n), $$
 with the initial conditions $f_0(1)=1$, $f_k(n)=0$ for $k>\lfloor n/2
 \rfloor$. The proof now follows from routine generating function
 arguments. $\ \Box$

Another proof of Theorem~\ref{thm:simsun} was given by Maria Monks
(private communication, 2008). She gives a bijection between simsun
permutations in $\sn$ and the oriented increasing binary trees on the
vertex set $[n]$ discussed in Section~\ref{subsec:flip}. Simsun
permutations have an interesting connection with the $cd$-index of
$\sn$, discussed in Section~\ref{sec5}. For some further work on
simsun permutations, see Chow and Shiu \cite{c-s} and Deutsch and
Elizalde \cite{d-e}. 

\subsection{Orbits of chains of partitions} \label{ssec:orb}
A \emph{partition} $\pi$ of the set $[n]$ is a collection
$\{B_1,\dots,B_k\}$ of nonempty subsets of $[n]$ (called the
\emph{blocks} of $\pi$) such that $\bigcup B_i=[n]$ and $B_i\cap B_j=
\emptyset$ if $i\neq j$. Let $\Pi_n$ denote the set of all partitions
of $[n]$.  If $\pi,\sigma\in \Pi_n$, then we say that $\pi$ is a
\emph{refinement} of $\sigma$, denoted $\pi\leq \sigma$, if every
block of $\pi$ is contained in a block of $\sigma$. The relation
$\leq$ is a partial order, so $\Pi_n$ becomes a partially ordered set
(poset). Note that the symmetric group $\sn$ acts on $\Pi_n$ in an
obvious way, viz., if $B=\{a_1,\dots,a_j\}$ is a block of $\pi$ and
$w\in\sn$, then $w\cdot B:=\{w(a_1),\dots, w(a_j)\}$ is a block of
$w\cdot \pi$.

Let $\mathcal{M}(\Pi_n)$ denote the set of all maximal chains of
$\Pi_n$, i.e., all chains
  $$ \pi_0<\pi_1<\cdots< \pi_{n-1}, $$
so that for all $0\leq i\leq n-2$, $\pi_{i+1}$ is obtained from
$\pi_i$ by merging two blocks of $\pi_i$. Thus $\pi_i$ has exactly
$n-i$ blocks. In particular, $\pi_0$ is the partition into $n$
singleton blocks, and $\pi_{n-1}$ is the partition into one block
$[n]$. The action of $\sn$ on $\Pi_n$ induces an action on
$\mathcal{M}(\Pi_n)$. For instance, when $n=5$ a class of orbit
representatives is given by the five chains below. We write e.g.\
$12\hy 34\hy 5$ for the partition $\{\{1,2\},\{3,4\},\{5\}\}$, and we
omit the first and last element of each chain.
  $$ \begin{array}{lclclclcl}
    12\hy 3\hy 4\hy 5 & < & 123\hy 4\hy 5 & < & 1234\hy 5\\[.05in]
    12\hy 3\hy 4\hy 5 & < & 123\hy 4\hy 5 & < & 123\hy 45\\[.05in]
    12\hy 3\hy 4\hy 5 & < & 12\hy 34\hy 5 & < & 125\hy 34\\[.05in]
    12\hy 3\hy 4\hy 5 & < & 12\hy 34\hy 5 & < & 12\hy 345\\[.05in]
    12\hy 3\hy 4\hy 5 & < & 12\hy 34\hy 5 & < & 1234\hy 5 \end{array}
  $$   

\begin{theorem} \label{thm:ptnorb}
The number of orbits of the action of $\sn$ on $\mathcal{M}(\Pi_n)$ is
the Euler number $E_{n-1}$.
\end{theorem}

\begin{figure}
\centering
\centerline{\includegraphics[width=12cm]{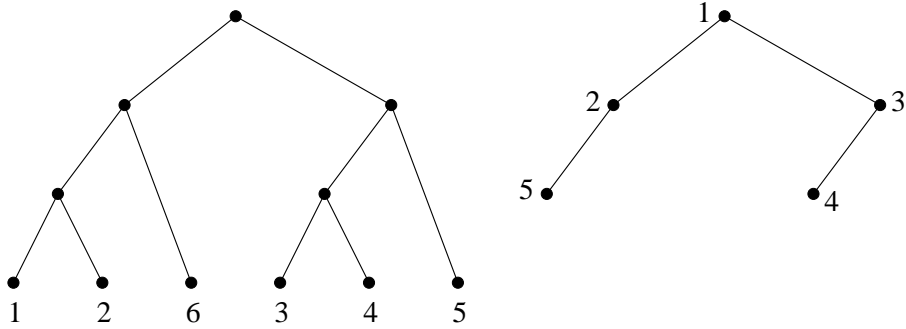}}
\caption{Two trees associated with a maximal chain in $\Pi_6$}
\label{fig:mtree}
\end{figure}

\proof
Given a maximal chain $\mathfrak{m}$ in $\mathcal{M}(\Pi_n)$, define a
binary tree with endpoints $1,2,\dots,n$ by the rule that for each
internal vertex $v$, at some point in the chain we merged together a
block consisting of the endpoints of the left subtree of $v$ with a
block consisting of the endpoints of the right subtree of $v$. For
instance, if the maximal chain is (omitting the first and last
elements) $12\hy 3\hy 4\hy 5\hy 6$, $12\hy 34\hy 5\hy 6$, $12\hy
345\hy 6$, $126\hy 345$, then the tree is given by
Figure~\ref{fig:mtree}(a). Label each internal vertex $v$ by $n-i$ if
that vertex was created at the $i$th step of the merging process and
delete the endpoints, as illustrated in Figure~\ref{fig:mtree}(b),
resulting in an increasing binary tree $T_{\mathfrak{m}}$ on vertices
$1,2,\dots,n$. The tree $T_{\mathfrak{m}}$ is well-defined up to flip
equivalence. Two maximal chains $\mathfrak{m}$ and $\mathfrak{m}'$
belong to the same $\sn$-orbit if and only $T_{\mathfrak{m}}$ and
$T_{\mathfrak{m}'}$ are flip equivalent, and the proof follows from
Theorem~\ref{prop:fostr}.
\qed

Theorem~\ref{thm:ptnorb} was first proved by Stanley
\cite[Thm.~7.7]{rs:aspects} by showing that the number of orbits
satisfied the recurrence \eqref{eq:eulerrec}. By elementary
representation theory, the number of orbits of $\sn$ acting on
$\mathcal{M}(\Pi_n)$ is the multiplicity of the trivial representation
in this action. This observation suggests the problem of decomposing
$\sn$-actions on various sets of chains in $\Pi_n$ into irreducible
representations. The first results in this direction appear in
\cite[{\S}7]{rs:aspects}. Many further results were obtained by
Sundaram \cite{sund}. Another such result is Theorem~\ref{thm:evenbl}
below. 

\subsection{Thickened zigzag tableaux} \label{subsec:zzt}
\begin{figure}
\centering
 \centerline{\includegraphics[width=6cm]{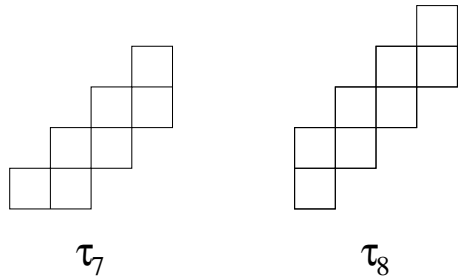}}
\caption{The zigzag shapes $\tau_7$ and $\tau_8$}
\label{fig:tau}
\end{figure}

For this subsection we assume familiarity with the theory of symmetric
functions such as developed in \cite{mac}\cite[Chap.~7]{ec2}. Let
$\tau_n$ be the border strip (or ribbon) corresponding to the
composition $\alpha=(1,2,2,\dots,2,j)$ of $n$, where $j=1$ if $n$ is
even and $j=2$ if $n$ is odd. Thus the (Young) diagram of $\tau_n$ has
a total of $n$ squares. Figure~\ref{fig:tau} shows the diagrams of
$\tau_7$ and $\tau_8$. We call $\tau_n$ a \emph{zigzag shape}.

\begin{figure}
\centering
 \centerline{\includegraphics[width=12cm]{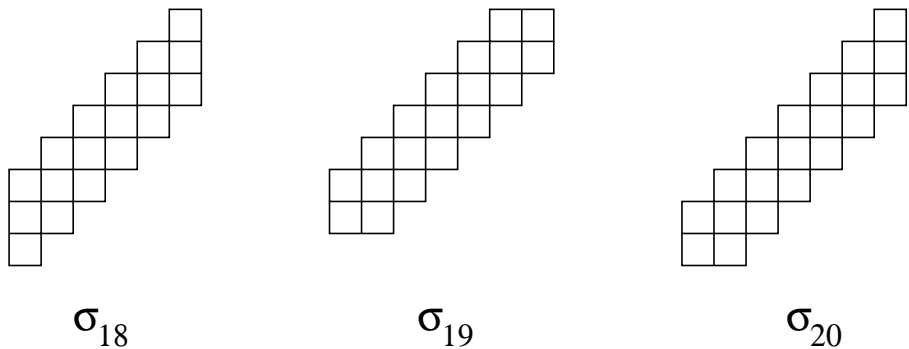}}
\caption{The thickened zigzag shapes $\sigma_{18}$, $\sigma_{19}$, and
              $\sigma_{20}$}
\label{fig:thick}
\end{figure}

Given any skew shape $\lambda/\mu$ of size $n$, let
$f^{\lambda/\mu}$ denote the number of standard Young tableaux (SYT)
of shape $\lambda/\mu$, i.e., the number of ways to put $1,2,\cdots,n$
into the squares of (the diagram of) $\lambda/\mu$, each number
$1,2,\dots,n$ occuring exactly once, so that the rows and colums are
increasing. If $T$ is an SYT of shape $\tau_n$, then reading the
numbers of $T$ from top-to-bottom and right-to-left gives a bijection
with reverse alternating permutations in $\sn$. Hence
  \beq f^{\tau_n} = E_n. \label{eq:taun} \eeq
Y. Baryshnikov and D. Romik \cite{b-r} give a surprising
generalization of equation~\eqref{eq:taun} in which the shapes
$\tau_n$ are ``thickened.'' We only mention the simplest case
here. The relevant shapes $\sigma_n$ (with $n$ squares) are
illustrated in Figure~\ref{fig:thick}, for each of the three case
$n\equiv 0,1,2\,(\mathrm{mod}\,3)$.

\begin{theorem} \label{thm:b-r}
We have
  \beas f^{\sigma_{3n-2}} & = & \frac{(3n-2)!E_{2n-1}}
          {(2n-1)!2^{2n-2}}\\
        f^{\sigma_{3n-1}} & = & \frac{(3n-1)!E_{2n-1}}
          {(2n-1)!2^{2n-1}}\\
       f^{\sigma_{3n}} & = & \frac{(3n)!(2^{2n-1}-1)E_{2n-1}}
          {(2n-1)!2^{2n-1}(2^{2n}-1)}. \eeas
\end{theorem}

The proof of Theorem~\ref{thm:b-r} generalizes a transfer operator
approach to alternating permutations due to Elkies \cite{elkies}. Is
there a bijective proof?

\subsection{M\"obius functions} \label{subsec:mu}
Let $P$ be a finite poset and Int$(P)$ the set of nonempty closed
invervals $[s,t]=\{u\st s\leq u\leq t\}$ of $P$. The \emph{M\"obius
  function} of $P$ (say over $\rr$) is the function $\mu\colon
\mathrm{Int}(P) \to\rr$ defined recursively as follows:
  \beas \mu(t,t) & = & 1,\ \ \mathrm{for\ all}\ t\in P\\
        \sum_{u\in[s,t]}\mu(s,u) & = & 0,\ \ \mathrm{for\ all}\ s<t\
         \mathrm{in}\ P. \eeas
Here we write $\mu(s,u)$ for $\mu([s,u])$. The M\"obius function has
many important properties and applications. See for instance
\cite[Chap.~3]{ec1} for more information. A number of posets have
Euler numbers as M\"obius function values (up to sign). We state the
two most significant such results here. Both these results are special
cases of much more general results that are related to topological
combinatorics and representation theory.

\begin{figure}
\centering
\centerline{\includegraphics[width=12cm]{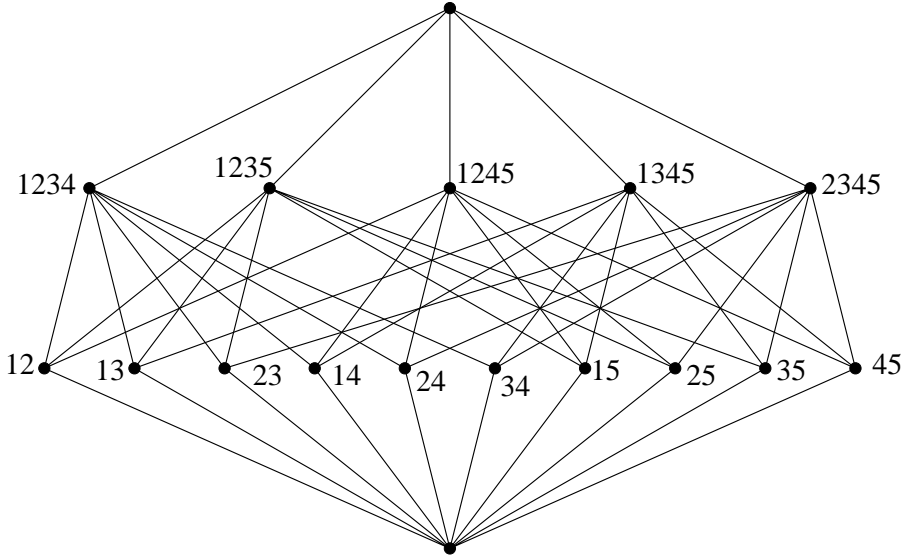}}
\caption{The poset $B_{5,2}$}
\label{fig:b52}
\end{figure}

For $n\geq 1$ let $B_{n,2}$ denote the poset of all subsets of
$[n]$ with an even number of elements, ordered by inclusion, with a
unique maximal element $\ho$ adjoined when $n$ is odd. Thus $B_{n,2}$
has a unique minimal element $\emptyset$ and a unique maximal element,
which we denote by $\ho$ even if $n$ is even. Figure~\ref{fig:b52}
shows $B_{5,2}$.   

\begin{theorem} \label{thm:mubne}
For the poset $B_{n,2}$ we have
   $$ \mu(\emptyset,\ho)=(-1)^{\lceil n/2\rceil} E_n. $$
\end{theorem}

Theorem~\ref{thm:mubne} is best understood in the general context of
rank-selected subposets of a graded poset. For an introduction to this
theory see \cite[{\S}3.12]{ec1}.

Now let $\Pi_{2n,2}$ denote the subposet of the poset $\Pi_{2n}$
defined in Section~\ref{ssec:orb} consisting of those partitions whose
block sizes are all even, with a minimal element $\hz$ adjoined. This
poset has a unique maximal element $\ho$ consisting of the partition
with one block $[2n]$.

\begin{theorem} \label{thm:evenbl}
For the poset $\Pi_{2n,2}$ we have
   $$ \mu(\hz,\ho) = (-1)^nE_{2n-1}. $$
\end{theorem}

Theorem~\ref{thm:evenbl} is due to G. Sylvester \cite{syl}. For some
generalizations see \cite{c-h-r} and \cite{rs:expst}.

\subsection{Polytope volumes}
Euler numbers occur as (normalized) volumes of certain convex
polytopes. The first polytope, which we call the \emph{zigzag
  order polytope} $\cp_n$, consists of all points
$x=(x_1,\dots,x_n)\in \rr^n$ satisfying
   $$ 0\leq x_i\leq 1,\ \ 1\leq i\leq n $$
   $$ x_1\geq x_2 \leq x_3 \geq \cdots x_n. $$
To compute its volume, for each alternating permutation
$w=a_1a_2\cdots a_n\in \sn$, let $w^{-1}=b_1b_2\cdots b_n$. Let 
   $$ \cp_w=\{ (x_1,\dots,x_n)\in\cp_n\st x_{b_1}\leq x_{b_2}\leq
       \cdots \leq x_{b_n}\}. $$
It is easy to see that each $\cp_w$ is a simplex with volume $1/n!$.
One can check using the theory of $P$-partitions \cite[{\S}4.5]{ec1}
that the $\cp_w$'s have disjoint interiors and union $\cp_n$. Hence
vol$(\cp_n)=E_n/n!$.

The second polytope is called the \emph{zigzag chain polytope}
$\cac_n$. It consists of all points $x=(x_1,\dots,x_n)\in \rr^n$
satisfying
   $$ x_i \geq 0,\ \ 1\leq i\leq n $$
   $$ x_i+x_{i+1}\leq 1,\ \ 1\leq i\leq n-1. $$  
The polytope $\cac_n$ first arose in \cite{rs:vol} and \cite{dob}.
It is also a special case of the \emph{Fibonacci polytopes} defined by
Rispoli \cite{rispoli}. A ``naive'' method for computing the volume is
the following. For $0\leq t\leq 1$ let
  \beq f_n(t) =\int_{x_1=0}^t\int_{x_2=0}^{1-x_1}\int_{x_3=0}^{1-x_2}
    \cdots \int_{x_n=0}^{1-x_{n-1}}dx_1\,dx_2\cdots dx_n. 
    \label{eq:int} \eeq
Clearly $f(1)=\mathrm{vol}(\cac_n)$. Differentiating
equation~\eqref{eq:int} yields $f'_n(t) = f_{n-1}(1-t)$. There are
various ways to solve this recurrence for $f_n(t)$ (with the initial
conditions $f_0(t)=1$ and $f_n(0)=0$ for $n>0$), yielding
   $$ \sum_{n\geq 0}f_n(t)x^n = (\sec x)(\cos(t-1)x+\sin tx). $$
Putting $t=1$ gives
   $$ \sum_{n\geq 0}f_n(1)x^n = \sec x+\tan x, $$
so we conclude that vol$(\cac_n)=E_n/n!$.

A more sophisticated proof uses the earlier obtained fact that
vol$(\cp_n)=E_n/n!$. Given $(x_1,\dots,x_n)\in\rr_n$, define
$\varphi(x_1,\dots,x_n)=(y_1,\dots,y_n)\in\rr^n$ by
  $$ y_i =\left\{ \begin{array}{rl} 1-x_i, & \mathrm{if}\ i\
      \mathrm{is\ odd}\\
     x_i, & \mathrm{if}\ i\ \mathrm{is\ even}. \end{array}
    \right. $$
It is easily checked that $\varphi$ is an affine transformation taking
$\cp_n$ onto $\cac_n$. Since the homogeneous part of $\varphi$ has
determinant $\pm 1$, it follows that $\varphi$ is a volume-preserving
bijection from $\cp_n$ onto $\cac_n$, so
vol$(\cac_n)=\mathrm{vol}(\cp_n) = E_n/n!$. This argument appeared in
Stanley \cite[Thm.~2.3 and Exam.~4.3]{rs:2pol}. Chebikin and
Ehrenborg\cite{cheb} compute the $f$-vector (which gives the number
of faces of each dimension) of a generalization of the polytopes
$\cp_n$. Since $\cp_n$ and $\cac_n$ are affinely equivalent, this
computation also gives the $f$-vector of $\cac_n$.

The polytope $\cac_n$ has an interesting connection to tridiagonal
matrices. An $n\times n$ matrix $M=(m_{ij})$ is \emph{tridiagonal} if
$m_{ij}=0$ whenever $|i-j|\geq 2$. Let $\mathcal{T}_n$ be the set of
all $n\times n$ tridiagonal doubly stochastic matrices $M$, i.e.,
$n\times n$ (real) tridiagonal matrices with nonnegative entries and
with row and column sums equal to 1.\ \ Thus $\mathcal{T}_n$ is a
convex polytope in a real vector space of dimension $n^2$ (or of
dimension $3n-2$ if we discard coordinates that are always 0). It is
easy to see that if we choose the $n-1$ entries
$m_{12}$, $m_{23}$, $\dots, m_{n-1,n}$ arbitrarily, then they
determine a unique tridiagonal matrix $M$ with row and column sums
1. Moreover, in order for $M$ to be doubly stochastic it is necessary
and sufficient that $m_{i,i+1}\geq 0$ and
  $$ m_{12}+m_{23}\leq 1,\ \ m_{23}+m_{34}\leq 1, \dots,\ 
       m_{n-2,n-1}+m_{n-1,n}\leq 1. $$
It follows that $\mathcal{T}_n$ is linearly equivalent to $\cac_{n-1}$
(in fact, $\mathcal{T}_n$ projects bijectively to $\cac_{n-1}$).
Moreover, the relative volume of $\mathcal{T}_n$ (volume normalized so
that a fundamental parallelopiped of the lattice
aff$(\mathcal{T}_n\cap \zz^{n^2}$ has volume 1, where aff denotes
affine span) is $E_{n-1}/(n-1)!$.

The $n\times n$ tridiagonal doubly-stochastic matrices form a face of
the \emph{Birkhoff polytope} $\mathcal{B}_n$ of all $n\times n$ doubly
stochastic matrices. Another face of $\mathcal{B}_n$ with an
interesting volume is the \emph{Chan-Robbins polytope}
$\mathcal{CR}_n$ defined by
  $$ \mathcal{CR}_n =\{ M=(m_{ij})\in\mathcal{B}_n\st m_{ij}=0\
           \mathrm{if}\ i-j\geq 2 \}, $$
for which vol$(\mathcal{CR}_n)=C_1 C_2\cdots C_n/\binom n2!$ (where
$C_i$ is a Catalan number) \cite{ba-ve}\cite{zeil}. The nice formulas
for the volumes of $\mathcal{T}_n$ and $\mathcal{CR}_n$ suggest the
problem of finding an interpolation between the two. For instance, for
$1\leq k\leq n-1$, can one compute the volume of the polytope
  $$ \mathcal{T}_{n,k} =\{ M=(m_{ij})\in\mathcal{B}_n\st m_{ij}=0\
           \mathrm{if}\ i-j\geq 2\ \mathrm{or}\ j-i>k \}? $$
Note that $\mathcal{T}_{n,n-1}=\mathcal{B}_n$ and $\mathcal{T}_{n,1}=
\mathcal{CR}_n$. 

\subsection{Singularities}
V. I. Arnold \cite{arnold} (see also \cite{arnold2} for a followup)
has shown that the Euler number $E_{n+1}$ is equal to the number of
components of the space of real polynomials $f(x) = x^n+a_1x^{n-1}
+\cdots + a_{n-1}x$ whose critical points (zeros of $f'(x)$) are all
real and whose $n-1$ critical values (the numbers $f(c)$ where $c$ is a
critical point) are all different. For instance, when $n=3$ the
polynomials $x^3+ax^2+bx$ form a real plane. The critical points are
real if and only if $b\leq a^2/3$. Two critical values coincide in the
region $b<a^2/3$ if and only if $b=a^2/4$ or $b=0$. These two curves
cut the region $b<a^2/3$ into $E_4=5$ components. Arnold interprets
this result in terms of morsifications of the function $x^{n+1}$; see
his paper for further details. Arnold goes on to deduce a number of
interesting properties of Euler numbers. He also extends the theory to
morsifications of the functions $x^n+y^2$ and $xy+y^n$, thereby
producing $B_n$ and $C_n$ analogues of Euler numbers (which correspond
to the root system $A_n$).

\section{Longest alternating subsequences} \label{sec3} 
Much work has been devoted to the properties of the length is$(w)$ of
the longest increasing subsequence of a permutation $a_1\cdots a_n$,
i.e., the largest $k$ for which there exist $i_1<\cdots<i_k$ and
$a_{i_1}<\cdots <a_{i_k}$. For a survey of this subject, see
\cite{rs:ids}. Two of the highlights of this subject are the
following. Let $E(n)$ denote the expected length 
of the longest increasing subsequence of $w\in\sn$ (with respect to
the uniform distribution on $\sn$). In symbols,
  $$ E(n) = \frac{1}{n!}\sum_{w\in\sn} \is(w). $$
It was shown by Vershik and Kerov \cite{v-k}, and the difficult part
of the argument independently by Logan and Shepp \cite{l-s}, that
  \beq E(n) \sim 2\sqrt{n}. \label{eq:en} \eeq
The notation $f(n)\sim g(n)$ means that $\lim_{n\to\infty}f(n)/g(n) =
1$. 

A far-reaching improvement of equation~\eqref{eq:en} was given by
Baik, Deift, and Johansson \cite{b-d-j}, namely, they determined the
(suitably scaled) limiting distribution of $\is(w)$ for $w\in\sn$ as
$n\to \infty$. Let $F(t)$ denote the \emph{Tracy-Widom distribution},
a probability distribution on $\rr$ first arising in the work of Tracy
and Widom on eigenvalues of random hermitian matrices \cite{t-w}. (We
will not define this distribution here.) Write $\is_n(w)$ for
$\is(w)$, where $w\in\sn$. The result of Baik, Deift, and Johansson
asserts that 
    \beq \lim_{n\to \infty} \mathrm{Prob}\left(
    \frac{\mathrm{is}_n(w)-2\sqrt{n}} 
     {n^{1/6}}\leq t\right)  =  F(t). \label{eq:bdj} \eeq
Here for each $n$ we are choosing $w\in\sn$ independently and
uniformly. 

We can ask whether similar results hold for \emph{alternating}
subsequences of $w\in\sn$. In particular, for $w\in\sn$ define
$\as(w)$ (or $\as_n(w)$ to make it explicit that $w\in\sn$) to be the
length of the longest alternating subsequence of $w$. For instance, if
$w=56218347$ then $\as(w)=5$, one alternating subsequence of longest
length being 52834. Our source for material in this section is
the paper \cite{rs:as}. 

It turns out that the behavior of $\as(w)$ is much simpler than that
of $\is(w)$. The primary reason for this is the following lemma, whose
straightforward proof we omit.

\begin{lemma} \label{lemma:key}
Let $w\in\sn$. Then there is an alternating subsequence of $w$ of
maximum length that contains $n$.
\end{lemma}

Lemma~\ref{lemma:key} allows us to obtain explicit formulas by
induction. More specifically, define   
  \beas a_k(n) & = & \#\{w\in\sn\st \as(w)=k\}\\[1em]
    b_k(n) & = & a_1(n)+a_2(n)+\cdots+a_k(n)\\
     & = & \#\{w\in\sn\st\as(w)\leq k\}. \eeas
For instance, $b_1(n)=1$, corresponding to the permutation
$1,2,\dots,n$, while $b_2(n)=2^{n-1}$, corresponding to the permutations
$u_1,u_2\dots,u_i,n,v_1,v_2,\dots,v_{n-i-1}$, where
$u_1<u_2<\cdots<u_i$ and $v_1>v_2>\cdots>v_{n-i-1}$. Using
Lemma~\ref{lemma:key}, we can obtain the following recurrence for the
numbers $a_k(n)$, together with the intial condition $a_0(0)=1$:
    \beq a_k(n+1) = \sum_{j=0}^{n} \binom nj
   \sum_{\substack{2r+s=k-1\\ r,s\geq
   0}}(a_{2r}(j)+a_{2r+1}(j))a_s(n-j).  
   \label{eq:aknrec} \eeq

This recurrence can be used to obtain the following generating
function for the numbers $a_k(n)$ and $b_k(n)$. No analogous formula
is known for increasing subsequences. 

\begin{theorem} \label{thm:main}
Let 
  \beas A(x,t) & = & \sum_{k,n\geq 0}a_k(n)t^k\frac{x^n}{n!}\\
   B(x,t) & = & \sum_{k,n\geq 0}b_k(n)t^k\frac{x^n}{n!}. \eeas
Set $\rho=\sqrt{1-t^2}$. Then
  \beas B(x,t) & = & \frac{ 2/\rho}{\displaystyle
      1-\frac{1-\rho}t e^{\rho 
      x}} -\frac{1}{\rho}\\[1em]
      A(x,t) & = & (1-t)B(x,t). \eeas
\end{theorem}

Many consequences can be derived from Theorem~\ref{thm:main}. In
particular, there are explicit formulas for $a_k(n)$ and $b_k(n)$.

\begin{cor} \label{cor:aknbkn}
For all $k,n\geq 1$ we have
  \begin{align} 
  b_k(n) & = \frac{1}{2^{k-1}} \sum_{\substack{r+2s\leq k\\ r\equiv
     k\,(\mathrm{mod}\,2)}} 
     (-2)^s \binom{k-s}{(k+r)/2}\binom ns r^n \label{eq:bkns}\\
   a_k(n) & = b_k(n)-b_{k-1}(n). \label{eq:akns} \end{align}
\end{cor}
 
For $k\leq 6$ we have
  \begin{align*}  b_2(n) & = 2^{n-1}\\
        b_3(n) & = \frac 14(3^n-2n+3)\\
    b_4(n) & = \frac 18(4^n-2(n-2)2^n)\\ 
    b_5(n) & = \frac{1}{16}(5^n-(2n-5)3^n+2(n^2-5n+5)) \\
    b_6(n) & = \frac{1}{32}(6^n-2(n-3)4^n+(2n^2-12n+15)2^n). 
\end{align*}

We can also obtain explicit formulas for the moments of $\as(w)$. For
instance, to obtain the mean (expectation)
  $$ D(n) = \frac{1}{n!}\sum_{w\in\sn} \as(w), $$
we compute
  \beas \sum_{n\geq 1}D(n)x^n & = & 
      \frac{\partial}{\partial t}A(x,1)\\ & = &
     \frac{6x-3x^2+x^3}{6(1-x)^2}\\ & = &
     x+\sum_{n\geq 2}\frac{4n+1}{6}x^n. \eeas
Thus 
   $$ D(n) = \frac{4n+1}{6},\ \ n\geq 2, $$
a remarkably simple formula. Note that (not surprisingly) $D(n)$ is
much larger than the expectation of is$(w)$, viz., $E(n)\sim
2\sqrt{n}$. Similarly one can obtain explicit
formulas for all the higher moments. In particular, the variance
  $$ V(n)
  =\frac{1}{n!}\sum_{w\in\sn}\left(\as(w)-D(n)\right)^2 $$
is given by 
   $$ V(n) =\frac{8}{45}n-\frac{13}{180},\ \ n\geq 4. $$
\indent Now that we have computed the mean and variance of $\as(w)$,
we can ask whether there is an ``alternating analogue'' of the
Baik-Deift-Johansson formula \eqref{eq:bdj}. In other words, can we
determine the scaled limiting distribution
   $$ K(t) = \lim_{n\to\infty} \mathrm{Prob}\left(
   \frac{\as_n(w)-2n/3}{\sqrt{n}} \leq t\right), $$
for $t\in \rr$? It turns out that the limiting distribution is
Gaussian. It is a consequence of results of Pemantle and Wilson
\cite{p-w} and Wilf \cite{wilf}, and was proved directly by Widom
\cite{widom}. More precisely, we have
  $$  K(t) = \frac{1}{\sqrt{\pi}} \int_{-\infty}^{t\sqrt{45}/4}
        e^{-s^2}ds. $$
Let us mention an observation of B\'ona that the statistic $\as(w)$ is
closely related to another statistic on permutations. Namely, an
\emph{alternating run} of a permutation $w\in\sn$ is a maximal factor
(subsequence of consecutive elements) that is increasing or 
decreasing. For instance, the permutation 64283157 has four
alternating runs, viz., 642, 28, 831, and 157. Let $g_k(n)$ be the
number of permutations $w\in\sn$ with $k$ alternating runs. Then
B\'ona's observation \cite{bona2} is that
  \beq a_k(n) = \frac 12(g_{k-1}(n)+g_k(n)),\ \ n\geq 2. 
    \label{eq:ag} \eeq
Hence all our results on as$(w)$ can be interpreted in terms of the
number of alternating runs of $w$. For some references to work on
alternating runs, see \cite[{\S}4]{rs:as}.

A comparison of results on is$(w)$ and $\as(w)$ suggests that it might
be interesting to interpolate between them. One possibility is the
following. Given $k\geq 1$, define a sequence $a_1a_2\cdots a_k$ of
integers to be \emph{$k$-alternating} if 
  $$ a_i>a_{i+1} \Leftrightarrow i\equiv 1\,(\mathrm{mod}\,k). $$
For instance, 61482573 is 3-alternating. A sequence is 2-alternating
if and only if it is alternating, and a sequence of length $n$ is
increasing if and only if it is $k$-alternating for some (or any)
$k\geq n-1$. What is the expected value $E_k(n)$ and limiting
distribution of the length of the longest $k$-alternating subsequence
of a random permutation $w\in\sn$? If $k$ is constant then most likely
$E_k(n) \sim c_k n$ for some constant $c_k$, while the limiting
distribution remains Gaussian (as for $k=2$). But what if $k$ grows
with $n$, e.g., $k=\lfloor \sqrt{n}\rfloor$? Is there a sharp cutoff
between the behavior $E_k(n)\sim cn$ and $E_k(n)\sim c\sqrt{n}$, or
is there a wide range of intermediate values? Similarly, do we get
limiting distributions other than Gaussian or Tracy-Widom? The same
questions can be asked if we replace $k$-alternating with the
condition that $a_i>a_{i+1}$ if and only if $\lfloor i/k\rfloor$ is
even, i.e., the permutation begins with $k-1$ descents, then $k-1$
ascents, etc. 

\section{Umbral enumeration of classes of alternating permutations}
\label{sec4} 
In this section we consider the enumeration of alternating
permutations having additional properties, such as having alternating
inverses or having no fixed points. The main tool will be a certain
character $\chi^{\tau_n}$ of the symmetric group $\sn$, first
considered by H. O. Foulkes \cite{foulkes2}\cite{foulkes}, whose
dimension in $E_n$. For the definition of $\chi^{\tau_n}$ we assume
familiarity with the theory of symmetric functions as in
Section~\ref{subsec:zzt}. 

For any skew shape $\lambda/\mu$ of size $n$ we can associate a
character $\chi^{\lambda/\mu}$ of $\sn$, e.g., by letting the
Frobenius characteristic ch$(\chi^{\lambda/\mu})$ be the skew Schur
function $s_{\lambda/\mu}$. In particular, taking $\tau_n$ to be the
zigzag shape of Section~\ref{subsec:zzt} gives the character
$\chi^{\tau_n}$. For any character $\chi$ of $\sn$ and
partition $\mu$ of $n$ we write $\chi(\mu)$ for $\chi(w)$, where $w$
is a permutation of cycle type $\mu$. The main result
\cite[Thm.~6.1]{foulkes}\cite[Exer.~7.64]{ec2} of Foulkes on the
connection between alternating permutations and representation theory
is the following.

\begin{theorem} \label{thm:foulkes}
 \emph{(a)} Let $\mu\vdash n$, where $n=2k+1$. Then
  $$ \chi^{\tau_n}(\mu) = \chi^{\tau'_n}(\mu) = \left\{ 
    \begin{array}{rl} 0, & \mbox{if $\mu$ has an even part}\\[.05in]
         (-1)^{k+r}E_{2r+1}, & \mbox{if $\mu$ has $2r+1$ odd parts
         and}\\  & \mbox{\ \ no even parts}. \end{array} \right. $$
 \emph{(b)} Let $\mu\vdash n$, where $n=2k$. Suppose that $\mu$ has
$2r$ odd parts and $e$ even parts. Then
   \beas \chi^{\tau_n}(\mu) & = &  (-1)^{k+r+e}E_{2r}\\[.05in]
     \chi^{\tau'_n}(\mu) & = &  (-1)^{k+r}E_{2r}. \eeas
\end{theorem}

Foulkes' result leads immediately to the main tool of this section.
We will use \emph{umbral notation} \cite{r-t}
for Euler numbers. In other words, any polynomial in $E$ is to be
expanded in terms of powers of $E$, and then $E^k$ is replaced by
$E_k$. The replacement of $E^k$ by $E_k$ is always the \emph{last} step in
the evaluation of an umbral expression.  For instance,
  $$ (E^2-1)^2 = E^4 -2E^2+1=E_4-2E_2+1=5-2\cdot 1+1=4. $$
Similarly,
  \beas (1+t)^E & = & 1+Et+\binom E2t^2+\binom E3 t^3+\cdots\\
    & = & 1+Et+\frac 12(E^2-E)t^2+\frac 16(E^3-3E^2+2E)t^3+\cdots\\
    & = & 1+E_1t+\frac 12(E_2-E_1)t^2+\frac 16(E_3-3E_2+2E_1)t^3+
       \cdots\\ & = & 
     1+1\cdot t+\frac 12(1-1)t^2+\frac 16(2-3\cdot 1+2\cdot 1)t^3
        + \cdots\\ & = & 1+t+\frac 16 t^3+\cdots. \eeas
If $f=f(x_1,x_2,\dots)$ is a
symmetric function then we use the notation $f[p_1,p_2,\dots]$ for $f$
regarded as a polynomial in the power sums. For instance, if
$f=e_2=\sum_{i<j}x_i x_j=\frac 12(p_1^2-p_2)$ then
  $$ e_2[E,-E,\dots]= \frac 12(E^2+E)=1. $$ 
We also write $\langle f,g\rangle$ for the standard (Hall) scalar
product of the symmetric functions $f$ and $g$.

\begin{theorem} \label{thm:main2}
Let $f$ be a homogenous symmetric function of degree $n$. If $n$ is
  odd then 
  \beq \langle f,s_{\tau_n}\rangle =\langle f,s_{\tau'_n}\rangle =
         f[E,0,-E,0,E,0,-E,\dots] \label{eq:feo} \eeq
If $n$ is even then 
 \beas \langle f,s_{\tau_n}\rangle & = & 
         f[E,-1,-E,1,E,-1,-E,1,\dots]\\
   \langle f,s_{\tau'_n}\rangle & = & 
         f[E,1,-E,-1,E,1,-E,-1,\dots]. \eeas
\end{theorem}

There are numerous results in the theory of symmetric functions that
express the number of permutations with certain properties as scalar
products of symmetric functions. Using Theorem~\ref{thm:main2} we
can obtain umbral generating functions for the number of alternating
permutations with various properties. Sometimes it is possible to
``deumbralize'' the generating function by expanding it explicitly in
powers of $E$ and then replacing $E^k$ with $E_k$. 

As an example of the procedure described above, let $f(n)$ denote the
number of permutations $w\in\sn$ such that both $w$ and $w^{-1}$ are
alternating. Such permutations are called \emph{doubly
  alternating}. As a special case of a well-known result of Foulkes 
\cite[Thm.~6.2]{foulkes2} (see also \cite[Cor.~7.23.8]{ec2}), we have
   $$ f(n) = \langle s_{\tau_n},s_{\tau_n}\rangle. $$
We can now use Theorem~\ref{thm:main2} to obtain (see \cite{rs:umbral}
for details)
    \beas \sum_{k\geq 0}f(2k+1)t^{2k+1} & = & \sum_{r\geq 0}E_{2r+1}^2
    \frac{L(t)^{2r+1}}{(2r+1)!}\\[.1in]
   \sum_{k\geq 0}f(2k)t^{2k} & = & \frac{1}{\sqrt{1-t^2}} \sum_{r\geq
     0}E_{2r}^2 \frac{L(t)^{2r}}{(2r)!}, \eeas
where 
      $$ L(t) = \frac 12\log\frac{1+t}{1-t} =
    t+\frac{t^3}{3}+\frac{t^5}{5}+\cdots. $$

A further class of alternating permutations that can be enumerated by
the above technique are those of a given cycle type, i.e., with
specified cycle lengths in their representation as a product of
disjoint cycles. We will only mention two special cases here. The
first case is when all cycles have length two. Such permutations are
just the fixed-point-free involutions. They necessarily have even
length. Thus let $g(n)$ (respectively, $g^*(n)$) denote the number of
fixed-point-free alternating involutions (respectively,
fixed-point-free reverse alternating involutions) in $\fs_{2n}$. Set
  $$ \begin{array}{lllll} G(t) & = & \sum_{n\geq 0} g(n)x^n & = & 1 +
    t +t^2 + 2t^3 + 5t^4 + 17t^5 + 72t^6+\cdots\\[.5em]
   G^*(t) & = & \sum_{n\geq 0} g^*(n)x^n & = & 1 + t^2
    + t^3 + 4t^4 + 13t^5 + 59t^6 +\cdots. \end{array} $$

\begin{theorem} \label{thm:mevenfm}
We have the umbral generating functions
  \beas G(t)  & = &  \left(\frac{1+t}{1-t}
    \right)^{(E^2+1)/4}\\[.05in] 
    G^\ast(t) & = & \displaystyle \frac{G(t)}{1+t}. \eeas
\end{theorem}

Ramanujan asserts in Entry 16 of his second notebook (see
\cite[p.~545]{berndt}) that as $t$ tends to $0+$,
  \beq 2\sum_{n\geq 0}(-1)^n\left( \frac{1-t}{1+t}\right)^{n(n+1)}
    \sim 1+t+t^2+2t^3+5t^4+17t^5+\cdots. \label{eq:berndt} \eeq
Berndt \cite[(16.6)]{berndt} obtains a formula for the complete
asymptotic expansion of $2\sum_{n\geq 0}(-1)^n\left(
\frac{1-t}{1+t}\right)^{n(n+1)}$ as $t\rightarrow 0+$. It is easy to
see that Berndt's formula can be written as $\left(\frac{1+t}{1-t}
\right)^{(E^2+1)/4}$ and is thus equal to
$G(t)$. Theorem~\ref{thm:mevenfm} therefore answers a question of
Galway \cite[p.~111]{galway}, who asks for a combinatorial
interpretation of the coefficients in Ramanujan's asymptotic
expansion. Is there a direct way to see the connection between
equation~\eqref{eq:berndt} and fixed-point-free alternating
involutions? Can equation~\eqref{eq:berndt} be generalized to involve
other classes of permutations?

\textsc{Note.} The following formula for $G(t)$ follows from
equation (\ref{eq:berndt}) and an identity of Ramanujan proved by
Andrews \cite[(6.3)$_\mathrm{R}$]{andrews}: 
  $$ F_2(t) = 2\sum_{n\geq 0}q^n\frac{\prod_{j=1}^n(1-q^{2j-1})}
    {\prod_{j=1}^{2n+1}(1+q^j)}, $$
where $q=\left(\frac{1-t}{1+t}\right)^{2/3}$. It is not hard to see
that this is a \emph{formal} identity, unlike the asymptotic identity
(\ref{eq:berndt}). 

The second class of alternating permutations with a given cycle type
that we will discuss are those that are single cycles. Let $b(n)$
(respectively, $b^*(n)$) denote the number of alternating
(respectively, reverse alternating) $n$-cycles in $\sn$. 

\begin{theorem} \label{thm:bn}
\emph{(a)} If $n$ is odd then
  $$ b(n) = b^\ast(n) = \frac 1n \sum_{d\mid
    n}\mu(d)(-1)^{(d-1)/2}E_{n/d}. $$ 
\emph{(b)} If $n=2^km$ where $k\geq 1$, $m$ is odd, and $m\geq 3$,
then
  $$ b(n) = b^*(n) = \frac 1n \sum_{d\mid m}\mu(d)E_{n/d}. $$
\emph{(c)} If $n=2^k$ and $k\geq 2$ then
  \beq b(n) = b^\ast(n) = \frac 1n(E_n-1). \label{eq:bnc} \eeq
\emph{(d)} Finally, $b(2)=1$, $b^\ast(2)=0$. 
\end{theorem}

Note that curious fact that $b(n)=b^*(n)$ except for $n=2$. Can
Theorem~\ref{thm:bn} be proved combinatorially, especially in the
cases when $n$ is a prime power (when the sums have only two terms)? 

It is immediate from Theorem~\ref{thm:bn} that $\lim_{n\to\infty} n
b(n)/E_n=1$. Thus as $n\to\infty$, a fraction $1/n$ of the alternating
permutations are $n$-cycles. Compare this with the simple fact that
(exactly) $1/n$ of the permutations $w\in\sn$ are $n$-cycles. We can
say that the properties of being an alternating permutation and an
$n$-cycle are ``asymptotically independent.'' What other class of
permutations are asymptotically independent from the alternating
permutations?

Closely related to the cycle types of alternating permutations is the
enumeration of alternating permutations with a given number of fixed
points, a question first raised by P. Diaconis. Thus let $d_k(n)$
(respectively, $d^*_k(n)$) denote the number of alternating
(respectively, reverse alternating) permutations in $\sn$ with exactly
$k$ fixed points. For a power series $F(t)=\sum a_nt^n$ we write
  $$ \begin{array}{lllll} \mathcal{O}_tF(t) & = & \frac 12(F(t)-F(-t))
        & = & \sum a_{2n+1}t^{2n+1}\\[1em]
    \mathcal{E}_tF(t) & = & \frac 12(F(t)+F(-t))
        & = & \sum a_{2n}t^{2n}. \end{array} $$

\begin{theorem} \label{prop:fixed}
We have 
  \beas \sum_{k,n\geq 0} d_k(2n+1)q^k t^{2n+1} & = & \mathcal{O}_t
   \frac{\exp(E(\tan^{-1}qt-\tan^{-1}t))}{1-Et}\\
  d^*_k(2n+1) & = & d_k(2n+1) \\
  \sum_{k,n\geq 0} d_k(2n)q^k t^{2n} & = & \mathcal{E}_t
  \sqrt{\frac{1+t^2}{1+q^2t^2}}
  \frac{\exp(E(\tan^{-1}qt-\tan^{-1}t))}{1-Et}\\ 
   \sum_{k,n\geq 0} d_k^*(2n)q^k t^{2n} & = &  \mathcal{E}_t
  \sqrt{\frac{1+q^2t^2}{1+t^2}}
  \frac{\exp(E(\tan^{-1}qt-\tan^{-1}t))}{1-Et}.
  \eeas
Equivalently, we have the non-umbral formulas
  \beas \sum_{k,n\geq 0} d_k(2n+1)q^k t^{2n+1} & = & 
    \sum_{\substack{i,j\geq 0\\ i\not\equiv j\,(\mathrm{mod}\,2)}}
     \frac{E_{i+j}}{j!}t^i(\tan^{-1}qt-\tan^{-1}t)^j\\
   \sum_{k,n\geq 0} d_k(2n)q^k t^{2n} & = & 
    \sqrt{\frac{1+t^2}{1+q^2t^2}}
    \sum_{\substack{i,j\geq 0\\ i\equiv j\,(\mathrm{mod}\,2)}}
     \frac{E_{i+j}}{j!}t^i(\tan^{-1}qt-\tan^{-1}t)^j\\
    \sum_{k,n\geq 0} d_k^*(2n)q^k t^{2n} & = & 
    \sqrt{\frac{1+q^2t^2}{1+t^2}}
    \sum_{\substack{i,j\geq 0\\ i\equiv j\,(\mathrm{mod}\,2)}}
     \frac{E_{i+j}}{j!}t^i(\tan^{-1}qt-\tan^{-1}t)^j.
  \eeas
\end{theorem}

Theorem~\ref{prop:fixed} may look a little daunting, but it can be
used to obtain some interesting concrete results. For instance, it is
not difficult to deduce from Theorem~\ref{prop:fixed} that
$d_0(n)=d_1(n)$ and $d_0^*(n)=d_1^*(n)$, yet another result that cries
out for a combinatorial proof. Moreover, we can obtain an asymptotic
expansion for the number of alternating derangements (permutations
without fixed points) in $\sn$, as follows.

\begin{cor} \label{cor:asy}
\rm{(a)} We have for $n$ odd the asymptotic expansion
  \beas d_0(n) & \sim & \frac 1e\left( E_n + a_1 E_{n-2} + a_2
  E_{n-4}+\cdots\right)\\ & = & \frac 1e\left(
  E_n+\frac 13 E_{n-2}- 
   \frac{13}{90}E_{n-4}+\frac{467}{5760}E_{n-6}+\cdots\right), 
    \eeas
where
   $$ \sum_{k\geq 0}a_k x^{2k} = \exp\left( 1-\frac 1x \tan^{-1}x
      \right). $$

\rm{(b)} We have for $n$ even the asymptotic expansion
 \beas d_0(n) & \sim & \frac 1e\left( E_n + b_1 E_{n-2} + b_2
  E_{n-4}+\cdots\right)\\ & = & \frac 1e\left(
  E_n+\frac 56 E_{n-2}- 
   \frac{37}{360}E_{n-4}+\frac{281}{9072}E_{n-6}+\cdots\right), 
    \eeas
where
   $$ \sum_{k\geq 0}b_k x^{2k} = \sqrt{1+x^2}\exp\left( 1-\frac 1x
      \tan^{-1}x \right). $$

    \rm{(c)} We have for $n$ even the asymptotic expansion
 \beas d_0^*(n) & \sim & \frac 1e\left( E_n + c_1 E_{n-2} + c_2
  E_{n-4}+\cdots\right)\\ & = & \frac 1e\left(
  E_n-\frac 16 E_{n-2}+ 
   \frac{23}{360}E_{n-4}-\frac{1493}{45360}E_{n-6}+\cdots\right), 
    \eeas
where
   $$ \sum_{k\geq 0}c_k x^{2k+1} = \frac{1}{\sqrt{1+x^2}}\exp\left(
      1-\frac 1x \tan^{-1}x \right). $$
\end{cor}

If $D(n)$ denotes the total number of derangements in $\sn$, then it
is a basic result from enumerative combinatorics that
$\lim_{n\to\infty} d(n)/n!=1/e$. On the other hand,
Corollary~\ref{cor:asy} shows that $\lim_{n\to\infty} d_0(n)/E_n =
\lim_{n\to\infty} d_0^*(n)/E_n =1/e$. In other words, the property of
being a derangement is asymptotically independent from being an
alternating permutation or from being a reverse alternating
permutation. 

A further consequence of Theorem~\ref{prop:fixed} concerns the
\emph{total} number of fixed points of alternating and reverse
alternating permutations in $\sn$. Again it is instructive to compare
with the case of all permutations. It is an immediate consequence of
Burnside's lemma (also known as the Cauchy-Frobenius lemma) and is
easy to see combinatorially that the total number of fixed points of
all permutations $w\in\sn$ is $n!$. For alternating permutations we
have the following result.

\begin{theorem} \label{thm:nofp}
Let $f(n)$ (respectively, $g(n)$) denote the total number of
fixed points of alternating (respectively, reverse
alternating) permutations in $\sn$. Then
   \beas  f(n) & = & \left\{ \begin{array}{rl}
     E_n-E_{n-2}+E_{n-4}-\cdots+(-1)^{(n-1)/2} E_1, & n\
   \mathrm{odd} \\[.05in]
    E_n-2E_{n-2}+2E_{n-4}-\cdots+(-1)^{(n-2)/2}2 E_2
  +(-1)^{n/2}, & n\ \mathrm{even}. \end{array} \right.\\ 
    g(n) & = & \left\{ \begin{array}{rl}
     E_n-E_{n-2}+E_{n-4}-\cdots+(-1)^{(n-1)/2} E_1, & n\ \mathrm{odd}
        \\[.05in]
    E_n-(-1)^{n/2}, & n\ \mathrm{even}. \end{array} \right. \eeas
\end{theorem}
Can these results be proved combinatorially?

There are umbral formulas for enumerating alternating and reverse
alternating permutations according to the length of their longest
increasing subsequence. These results are alternating analogues of the
following result of Gessel \cite[p.~280]{gessel}. Let $u_k(n)$ be the
number of permutations $w\in\sn$ whose longest increasing subsequence
has length at most $k$. For instance, $u_2(n)=C_n$, a Catalan number
\cite[Exer.~6.19(ee)]{ec2}. Define
   \beas U_k(x) & = & \sum_{n\geq 0} u_k(n)\frac{x^{2n}}{n!^2},\ k\geq
      1\\ 
    I_i(2x) & = & \sum_{n\geq 0}\frac{x^{2n+i}}{n!\,(n+i)!},\
       i\geq 0. \eeas 
The function $I_i$ is the \emph{hyperbolic Bessel function} of the
first kind of order $i$. Gessel then showed that
  \beq U_k(x) = \det\left( I_{|i-j|}(2x)\right)_{i,j=1}^k. 
        \label{eq:gesselis} \eeq
We state without proof the alternating analogue of
equation~\eqref{eq:gesselis}. It is proved by applying the umbral
techniques of \cite{rs:umbral} to a symmetric function identity
\cite[Thm.~16]{gessel} of Gessel. For any integer $i$ we use the
notation $D^iF(x)$ for the $n$th formal derivative of the power series
$F(x) = \sum_{n\geq 0}a_n x^n$, where in particular
  $$ D^{-1}F(x) = \sum_{n\geq 0}a_n\frac{x^{n+1}}{n+1}, $$
and $D^{-i-1}= D^{-1}D^{-i}$ for all $i\geq 1$.

\begin{theorem} \label{thm:altinc}
Let $v_k(n)$ (respectively, $v'_k(n)$) denote the number of
alternating (respectively, reverse alternating) permutations $w\in\sn$
whose longest increasing subsequence has length at most $n$. Let
$\exp\left(E\tan^{-1}(x)\right)=\sum_{n\geq 0}c_n(E)x^n$. Define
  \beas A_1(x) & = & \sum_{n\geq 0}c_n(E)\frac{x^n}{n!}\\
       A_2(x) & = & \sqrt{1+x^2}A_1(x)\\
       A_3(x) & = & \frac{A_1(x)}{\sqrt{1+x^2}}, \eeas
and for $1\leq r\leq 3$ and $k\geq 1$ define
       $$ B_{r,k}(x) = \det(D^{j-i}A_r(x))_{i,j=1}^k. $$
We then have:
  \begin{itemize} \item If $n$ is odd, then the coefficient of
    $x^n/n!$ in the umbral evaluation of $B_{1,k}(x)$ is $v_k(n) =
    v'_k(n)$.
  \item If $n$ is even, then the coefficient of
    $x^n/n!$ in the umbral evaluation of $B_{2,k}(x)$ is $v_k(n)$.      
  \item If $n$ is even, then the coefficient of
    $x^n/n!$ in the umbral evaluation of $B_{3,k}(x)$ is $v'_k(n)$. 
  \end{itemize}
\end{theorem}

Note that it follows from Theorem~\ref{thm:altinc} that
$v_k(2m+1)=v'_k(2m+1)$.  This result is also an immediate consequence
of the involution~\eqref{eq:inv} followed by reversal of the
permutation. 

It is known (see \cite{rs:catadd}) that $v_2(n)=C_{\lceil
  (n+1)/2\rceil}$ for $n\geq 3$ and $v'_2(n)=C_{\lceil n/2\rceil}$ for
$n\geq 1$, where $C_i=\frac{1}{i+1}\binom{2i}{i}$ (a Catalan number).
Can these formulas be deduced directly from the case $k=2$ of
Theorem~\ref{thm:altinc}? Similarly, J. Lewis \cite{lewis2} has shown
that  
  $$ v'_3(n) = \left\{ \begin{array}{rl}
     f^{(m,m,m)}, & n=2m\\
     f^{(m-1,m,m+1)}, & n=2m+1, \end{array} \right. $$
No such results are known for $v_4(n)$ or $v_3(n)$.

Is there an asymptotic formula for the expected length of the longest
increasing subsequence of an alternating permutation $w\in\sn$ as
$n\to\infty$, analogous to the result \eqref{eq:en} of Logan-Shepp and
Vershik-Kerov for arbitrary permutations?  Or more strongly, is there
a limiting distribution for the (suitably scaled) length of the
longest increasing subsequence of an alternating permutation $w\in\sn$
as $n\to\infty$, analogous to the Baik-Deift-Johansson theorem
\eqref{eq:bdj} for arbitrary permutations? Two closely related
problems are the following. Let $a^\lambda$ denote the number of SYT
of shape $\lambda\vdash n$ and descent set $\{1,3,5,\dots\}\cap
[n-1]$, as defined in \cite[p.~361]{ec2}. (These are just the SYT $Q$
that arise from alternating permutations $w\in\sn$ by applying the RSK
algorithm $w\mapsto (P,Q)$.) Let $f^\lambda$ denote the number of SYT
of shape $\lambda$, given by the famous \emph{hook length formula}
\cite[Cor.~7.21.6]{ec2}. What is the (suitably scaled) limiting shape
of $\lambda$ as $n\to\infty$ that maximizes $a^\lambda$ and similarly
that maximizes $a^\lambda f^\lambda$? (For the shape that maximizes
$f^\lambda$, see \cite[{\S}3]{rs:ids}. Could this shape also maximize
$a^\lambda$?) 

\textsc{Note.} Non-umbral methods have been used to enumerate certain
classes of (reverse) alternating permutations. A \emph{Baxter
  permutation} (originally called a 
\emph{reduced} Baxter permutation) is a permutation $w\in \sn$
satisfying: if $w(r)=i$ and $w(s)=i+1$, then there is a $k_i$ between
$r$ and $s$ (i.e., $r\leq k_i\leq s$ or $s\leq k_i\leq r$) such that
$w(t)\leq i$ if $t$ is between $r$ and $k_i$, while $w(t)\geq i+1$ if
$t$ is between $k_i+1$ and $s$. For instance, all permutations
$w\in\fs_4$ are Baxter permutations except 2413 and 3142.
Cori, Dulucq, and Viennot \cite{c-d-v} showed that the number of
reverse alternating Baxter permutations in $\fs_{2n}$ is $C_n^2$, and
in $\fs_{2n+1}$ is $C_nC_{n+1}$. Guibert and Linusson \cite{gu-li}
showed that the number of reverse alternating Baxter permutations
$\sn$ whose inverses are also reverse alternating is just $C_n$. A
different proof was later given by Min and Park \cite{mi-pa}.  Some
related work on enumerating special classes of alternating
permutations was done by Dulucq and Simion \cite{d-s}. A popular topic
in permutation enumeration is \emph{pattern avoidance}. See
\cite[Ch.~5]{bona} for an introduction. A special case deals with
permutations with no increasing subsequence of length $k$. Mansour
\cite{mansour} inaugurated the study of pattern avoidance in
alternating permutations. For some additional papers on this topic see
Hong \cite{hong}, Lewis \cite{lewis}\cite{lewis2}, and Ouchterlony
\cite{ouch}. 

\section{The $cd$-index of $\sn$} \label{sec5}
Let $w=a_1\cdots a_n\in\sn$. The \emph{descent set} $D(w)$ of $w$ is
defined by
  $$ D(w) = \{i\st a_i>a_{i+1}\}\subseteq [n-1]. $$
A permutation $w$ is thus alternating if $D(w)=\{1,3,5,\dots\}\cap
[n-1]$ and reverse alternating if $D(w)=\{2,4,6,\dots\}\cap[n-1]$. For
$S\subseteq [n-1]$ let 
  $$ \beta_n(S) =\#\{w\in\sn\st D(w)=S\}. $$
The numbers $\beta_n(S)$ are fundamental invariants of $\sn$ that
appear in a variety of combinatorial, algebraic, and geometric
contexts. In this section we explain how alternating permutations and
Euler numbers are related to the more general subject of permutations
with a fixed descent set.

We first define for fixed $n$ a noncommutative generating function for
the numbers $\beta_n(S)$. Given a set $S\subseteq [n-1]$, define its
\emph{characteristic monomial} (or \emph{variation}) to be the
noncommutative monomial  
   \beq u_S=e_1 e_2\cdots e_{n-1}, \label{eq:usvar} \eeq
where
    $$ e_i =\left\{ \begin{array}{rl}
     a, & \mathrm{if}\ i\not\in S\\
     b, & \mathrm{if}\ i\in S. \end{array} \right. $$
For instance, $D(37485216)=\{2,4,5,6\}$, so
$u_{D(37485216)}=ababbba$. Define 
  \bea \Psi_n = \Psi_n(a,b) & = & \sum_{w\in\sn} u_{D(w)} \nonumber\\
    & = & \sum_{S\subseteq [n-1]} \beta_n(S)u_S. \label{eq:psin} \eea
Thus $\Psi_n$ is a noncommutative generating function for the numbers
$\beta_n(S)$. For instance, 
  \beas \Psi_1 & = & 1\\ \Psi_2 & = & a+b\\
    \Psi_3 & = & a^2+2ab+2ba+b^2\\
   \Psi_4 & = & a^3+3a^2b+5aba+3ba^2+3ab^2+5bab+3b^2a+b^3. \eeas  
The polynomial $\Psi_n$ is called the \emph{$ab$-index} 
of the symmetric group $\sn$.

The main result of this section is the following.

\begin{theorem} \label{thm:cd}
There exists a polynomial $\Phi_n(c,d)$ in the noncommuting variables
$c$ and $d$ such that
  $$ \Psi_n(a,b) = \Phi_n(a+b,ab+ba). $$
\end{theorem}

The polynomial $\Phi_n(c,d)$ is called the \emph{$cd$-index} of
$\sn$. For instance, we have
  $$ \Psi_3(a,b) =a^2+2ab+2ba+b^2=(a+b)^2+(ab+ba), $$
so $\Phi_3(c,d)=c^2+d$. Some values of $\Phi_n(c,d)$ for small $n$ are
as follows:
   \beas \Phi_1 & = & 1\\ \Phi_2 & = & c\\
    \Phi_3 & = & c^2+d\\ \Phi_4 & = & c^3+2cd+2dc\\ 
    \Phi_5 & = & c^4+3c^2d+5cdc+3dc^2+4d^2\\
    \Phi_6 & = & c^5+4c^3d+9c^2dc+9cdc^2+4dc^3+12cd^2+10dcd+
       12d^2c. \eeas
If we define $\deg(c)=1$ and $\deg(d)=2$, then the number of
$cd$-monomials of degree $n-1$ is the Fibonacci number $F_n$. All
these monomials actually appear in $\Phi_n(c,d)$ (as in evident from
the discussion below). Thus $\Phi_n(c,d)$ has $F_n$ terms, compared
with $2^{n-1}$ terms for $\Psi_n(a,b)$.

There are several known proofs of Theorem~\ref{thm:cd}. Perhaps the
most natural approach is to define an equivalence relation on $\sn$
such that for each equivalence class $C$, we have that $\sum_{w\in C}
u_{D(w)}$ is a monomial in $c=a+b$ and $d=ab+ba$. Such a proof was
given by G.~Hetyei and E.~Reiner \cite{he-re}. We will simply define
the equivalence relation here. An exposition appears in Stanley
\cite[{\S}1.6]{ec1-2}. 

\begin{figure}
\centering
 \centerline{\includegraphics[width=10cm]{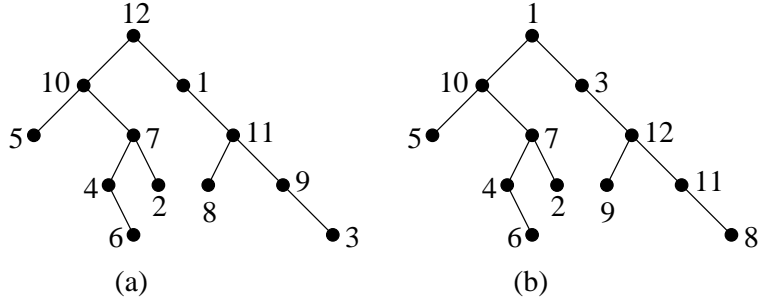}}
\caption{(a) The min-max tree
  $M=M(5,10,4,6,7,2,12,1,8,11,9,3)$; (b) 
  The transformed tree $\psi_7M=M(5,10,4,6,7,2,1,3,9,12,11,8)$}
\label{fig:mmtreex}
\end{figure}

We first define the \emph{min-max tree} $M(w)$ associated with a
sequence $w=a_1a_2\cdots a_n$ of distinct integers as follows. First,
$M(w)$ is a binary tree with vertices labelled $a_1, a_2,\dots,a_n$.
Let $j$ be the least integer for which \emph{either}
$a_j=\min\{a_1,\dots,a_n\}$ or $a_j=\max\{a_1,\dots,a_n\}$. Define
$a_j$ to be the root of $M(w)$.  Then define (recursively)
$M(a_1,\dots, a_{j-1})$ to be the left subtree of $a_j$, and
$M(a_{j+1},\dots,a_n)$ to be the right subtree.
Figure~\ref{fig:mmtreex}(a) shows $M(5,10,4,6,7,2,12,1,8,11,9,3)$.
Note that no vertex of a min-max tree $M(w)$ has only a left
successor. Note also that every vertex $v$ is either the minimum or
maximum element of the subtree with root $v$.

Given the min-max tree $M(w)$ where $w=a_1\cdots a_n$, we will define
operators $\psi_i$, $1\leq i\leq n$, that permute the labels of
$M(w)$, creating a new min-max tree $\psi_i M(w)$.  The operator
$\psi_i$ only permutes the label of the vertex of $M(w)$ labelled
$a_i$ and the labels of the right subtree of this vertex. (Note that
the vertex labelled $a_i$ depends only on $i$ and the tree $M(w)$, not
on the permutation $w$.) All other vertices are fixed by $\psi_i$. In
particular, if $a_i$ is an endpoint then $\psi_i M(w)=M(w)$.  Suppose
that $a_i$ is the minimum element of the subtree $M_{a_i}$ with root
$a_i$. Then replace $a_i$ with the \emph{largest} element of
$M_{a_i}$, and permute the remaining elements of $M_{a_i}$ so that
they keep their same relative order.  This defines $\psi_i M(w)$.
Similarly suppose that $a_i$ is the maximum element of the subtree
$M_{a_i}$ with root $a_i$. Then replace $a_i$ with the \emph{smallest}
element of $M_{a_i}$, and permute the remaining elements of $M_{a_i}$
so that they keep their same relative order. Again this defines
$\psi_i M(w)$.  Figure~\ref{fig:mmtreex}(b) shows that $\psi_7
M(5,10,4,6,7,2,12,1,8,11,9,3) = M(5,10,4,6,7,2,1,3,9,12,11,8)$. We
have $a_7=12$, so $\psi_7$ permutes vertex 12 and the vertices on the
right subtree of 12. Vertex 12 is replaced by 1, the smallest vertex
of the right subtree. The remaining elements $1,3,8,9,11$ get replaced
with $3,8,9,11,12$ in that order.

Let us call two permutations $v,w\in\sn$ \emph{equivalent}, denoted
$v\stackrel{M}{\sim} w$, if their min-max trees $T(v)$ and $T(w)$ can
be obtained from each other by applying a sequence of $\psi_i$'s.
Clearly $\stackrel{M}{\sim}$ is an equivalence relation on $\sn$.  Let
$c,d,e$ be noncommutative indeterminates, and let $w=a_1a_2\cdots
a_n\in\sn$.  For $1\leq i\leq n$ define
  $$ f_i =f_i(w)=\left\{ \begin{array}{rl}
    c, & \mbox{if\ $a_i$ has only a right child in}\ M(w)\\[.05in]
    d, & \mbox{if\ $a_i$ has a left and right child}\\[.05in]
    e, & \mbox{if\ $a_i$ is an endpoint.}
   \end{array} \right. $$
Let $\Phi'_w=\Phi'_w(c,d,e)=f_1f_2\cdots f_n$, and let
$\Phi_w=\Phi_w(c,d)= \Phi'(c,d,1)$, where 1 denotes the empty word. In
other words, $\Phi_w$ is obtained from $\Phi'_w$ by deleting the
$e$'s. For instance, consider the permutation
$w=5,10,4,6,7,2,12,1,8,11,9,3$  
of Figure~\ref{fig:mmtreex}. The degrees (number of children) of the
vertices $a_1,a_2,\dots,a_{12}$ are 0, 2, 1, 0, 2, 0, 2, 1, 0, 2, 1,
0, respectively. Hence
  \bea \Phi'_{w} & = & edcededcedce \nonumber\\
     \Phi_w & = & dcddcdc. \label{eq:phiwex} \eea
It is clear that if $v\stackrel{M}{\sim} w$, then $\Phi'_v=\Phi'_w$
and $\Phi_v=\Phi_w$, since $\Phi'_w$ depends only on $M(w)$ regarded
as an \emph{unlabelled} tree. 
   
The main result on min-max trees is the following.

\begin{theorem} \label{thm:minmax}
For any $w\in\sn$ we have 
   $$ \sum_{v\stackrel{M}{\sim} w}u_{D(v)} = \Phi_w(a+b,ab+ba). $$
\end{theorem}

Theorem~\ref{thm:minmax} shows that $\stackrel{M}{\sim}$ is precisely
the equivalence relation we asked for in order to prove the existence
of the $cd$-index $\Phi_n(c,d)$ (Theorem~\ref{thm:cd}). Thus not only
have we shown the existence of $\Phi_n(c,d)$, but also we have shown
that the coefficients are nonnegative. It is reasonable to ask whether
there is a more ``direct'' description of the coefficients. Such a
description was first given by D. Foata and M.-P.\ Sch\"utzenberger
\cite{fo-sc} in terms of the \emph{Andr\'e permutations} mentioned in
Section~\ref{ssec:simsun}. We state here the analogous result for
simsun permutations (as defined in Section~\ref{ssec:simsun}), due to
R. Simion and S. Sundaram.

\begin{theorem} \label{thm:cdcoef}
Let $\mu$ be a monomial of degree $n-1$ in the noncommuting variables
$c,d$, where $\deg(c)=1$ and $\deg(d)=2$. Replace each $c$ in $\mu$
with 0, each $d$ with 10, and remove the final 0. We get the
characteristic vector of a set $S_\mu\subseteq [n-2]$. Then the
coefficient of $\mu$ in $\Phi_n(c,d)$ is equal to the number of simsun
permutations in $\fs_{n-1}$ with descent set $S_\mu$.
\end{theorem}

For example,  if $\mu=cd^2c^2d$ then we get the characteristic vector
$01010001$ of the set $S_\mu=\{2,4,8\}$. Hence the coefficient of
$cd^2c^2d$ in $\Phi_{10}(c,d)$ is equal to the number of simsun
permutations in $\fs_9$ with descent set $\{2,4,8\}$.  

Note that every $cd$-monomial, when expanded in terms of
$ab$-monomials, is a sum of distinct monomials including
$bababa\cdots$ and $ababab\cdots$. These monomials correspond to
descent sets of alternating and reverse alternating permutations,
respectively. Hence $\Phi_n(1,1)=E_n$. This fact also follows
immediately from Theorems~\ref{thm:simsun} and \ref{thm:cdcoef}. 

Extending the reasoning of the previous paragraph gives a nice result
on inequalities among the numbers $\beta_n(S)$, originally due to
Niven \cite{niven} and de Bruijn \cite{bruijn} by different
methods. The proof given below first appeared in Stanley
\cite[Thm.~2.3(b)]{rs:cd} in a more general context.

 Given $S\subseteq[n-1]$, define
$\omega(S)\subseteq [n-2]$ by the condition $i\in\omega(S)$ if and
only if exactly one of $i$ and $i+1$ belongs to $S$, for $1\leq i\leq
n-2$. For instance, if $n=9$ and $S=\{2,4,5,8\}$, then
$\omega(S)=\{1,2,3,5,7\}$. Note that 
 \beq \omega(S)=[n-2]\ \Longleftrightarrow\
    S=\{1,3,5,\dots\}\cap[n-1]\ 
  \mathrm{or}\ S=\{2,4,6,\dots\}\cap[n-1]. \label{eq:omegamax} \eeq 

\begin{prop} \label{prop:betaineq}
Let $S,T\subseteq [n-1]$. If $\omega(S)\subset \omega(T)$, then
$\beta_n(S)<\beta_n(T)$. 
\end{prop}

\proof
Suppose that $\omega(S)\subseteq \omega(T)$. It is easy to check that
if $\mu$ is a $cd$-monomial such that the expansion of
$\mu(a+b,ab+ba)$ contains the term $u_T$ (necessarily with coefficient
1), then it also contains the term $u_S$. Since $\Phi_n(c,d)$ has
nonnegative coefficients, it follows that $\beta_n(S)\leq
\beta_n(T)$. 

Now assume that $S$ and $T$ are any subsets of $[n-1]$ for which
$\omega(S)\subset \omega(T)$ (strict containment). We can easily find
a $cd$-monomial $\Phi_w$ for which $\omega(T)\supseteq \omega (S_w)$ but
$\omega(S)\not\supseteq \omega(S_w)$. For instance, if $i\in
\omega(T)-\omega(S)$ then let $\Phi_w=c^{i-1}dc^{n-2-i}$, so
$\omega(S_w)=\{i\}$. It follows that $\beta_n(S) <\beta_n(T)$.  \qed

\begin{cor} \label{cor:betamax}
Let $S\subseteq [n-1]$. Then $\beta_n(S)\leq E_n$, with equality if
and only if $S=\{1,3,5,\dots\}\cap[n-1]$ or $S=\{2,4,6,\dots\}
\cap[n-1]$. 
\end{cor}

\proof
Immediate from Proposition~\ref{prop:betaineq} and
equation~\eqref{eq:omegamax}. 
\qed

Corollary~\ref{cor:betamax} can be rephrased as follows. If we pick a
permutation $w\in\sn$ at random (uniformly) and must predict its
descent set, then it is best to bet either that $w$ is alternating or
is reverse alternating. By equation~\eqref{eq:eulasymp} the probability
of success will be about $2(2/\pi)^{n+1}\approx
2(0.6366\cdots)^{n+1}$. 

An interesting generalization of Corollary~\ref{cor:betamax} is due to
K. Saito \cite[Thm.~3.2]{saito}. Let $T$ be a tree on an $n$-element
vertex set $V$. Let $\lambda$ be a labeling of $V$ with the numbers
$1,2,\dots,n$, and let $\fo_\lambda$ be the orientation of (the edges
of) $T$ obtained by orienting $u\rightarrow v$ if $uv$ is an edge of
$T$ with $\lambda(u)<\lambda(v)$. For each orientation $\fo$ of $T$,
let $\beta(\fo)$ be the number of labelings $\lambda$ for which $\fo =
\fo_\lambda$.

\begin{prop} \label{prop:saito}
The orientations $\fo$ that maximize $\beta(\fo)$ are the two
``bipartite orientations,'' i.e., those containing no directed path of
length 2. 
\end{prop}

Corollary~\ref{cor:betamax} is equivalent to the special case of
Proposition~\ref{prop:saito} for which $T$ is a path. The labelled
trees that produce bipartite orientations are known as
\emph{alternating } or \emph{intransitive} trees. They first appeared
in the work of Gelfand, Graev, and Postnikov \cite{g-g-p} in
connection with the theory of hypergeometric functions. They were
enumerated by Postnikov \cite{post} and are also connected with
counting regions of certain hyperplane arrangements
\cite[{\S}8.1]{po-rs}.

\textsc{Acknowledgment.} I am grateful to Joel Lewis for his careful
proofreading of this paper. 

\bibliographystyle{amsplain}

\end{document}